\documentclass[a4paper,11pt]{amsart}

\usepackage[english]{babel}
\usepackage[T1]{fontenc}
\usepackage{lmodern}

\usepackage[utf8]{inputenc}
\usepackage{amsmath,amssymb,mathrsfs}
\usepackage{amscd}
\usepackage{pgf,tikz}
\usepackage{hyperref}

\usepackage{a4wide}

\usepackage[utf8]{inputenc}
\usepackage{amsmath,amssymb,mathrsfs}
\usepackage{amscd}
\usepackage{pgf,tikz}
\usetikzlibrary{arrows}
\usetikzlibrary{patterns}
\numberwithin{equation}{section}
\newtheorem{theo}{Theorem}[section]
\newtheorem{cnj}{\textbf{\textup{Conjecture}}}[section]
\newtheorem{cor}{Corollary}[section]
\newtheorem{prop}{Proposition}[section]
\newtheorem{lem}{Lemma}[section]
\newtheorem{fe}{Definition}[section]
\newtheorem{rem}{Remark}[section]
\numberwithin{figure}{section}

\newcommand{\hn}{{P_n}}
\newcommand{\hd}{{P_2}}
\newcommand{\htr}{{P_3}}

\newcommand{\bn}{{\mathbb{N}}}
\newcommand{\bz}{{\mathbb{Z}}}
\newcommand{\br}{{\mathbb{R}}}
\newcommand{\bc}{{\mathbb{C}}}

\newcommand{\sqq}{{\mathcal{S}_q}}
\newcommand{\sd}{{\mathcal{S}_2}}

\newcommand{\sinf}{{\mathcal{S}_\infty}}
\newcommand{\one}{{\bf 1}}

\newcommand{\xp}{X_\perp}
\newcommand{\xpa}{X_\parallel}

\newenvironment{prof}
	{\textit{\textbf{Proof.}}}
	{\hfill $\square$\vskip 8pt}

\title[Distribution of complex eigenvalues]
{A simple criterion for the existence of nonreal eigenvalues
for a class of 2D and 3D Pauli operators}
\author{Diomba \textsc{Sambou}}

\address{Departamento de Matem\'aticas, Facultad de Matem\'aticas,
 Pontificia Universidad Cat\'olica de Chile, Vicu\~na Mackenna 4860, 
 Santiago de Chile}
 
\email{disambou@mat.uc.cl}

\thanks{The author is partially supported by the Chilean 
Program \textit{N\'ucleo Milenio de F\'isica Matem\'atica
RC$120002$}. The author gratefully acknowledges the many 
helpful suggestions of V. Bruneau during the preparation of 
the paper. The author should like to thank R. Novák for 
bringing to his attention the reference \cite{koc}.}

\begin{document}


\begin{abstract}
In this work, we investigate the discrete spectrum generated by complex 
matrix-valued perturbations for a class of 2D and 3D Pauli operators 
with nonconstant magnetic fields. We establish a simple criterion for 
the potentials to produce discrete spectrum near the low ground energy
of the operators. Moreover, in case of creation of nonreal eigenvalues, 
this criterion specifies also their location.



\end{abstract}

\maketitle


\section{Introduction}\label{s1}

\subsection{Description of the models}

We consider in this article $n$-dimensional Pauli operators 
$\hn(b,V)$, $n = 2$, $3$, defined as follows. Denote by 
$\xp := (x,y)$ the usual variables on $\br^2$ and by $X := 
(\xp,\xpa)$ those on $\br^3$. For $\textbf{x} = \xp \in \br^2$ 
or $\textbf{x} = X \in \br^3$, let
\begin{equation}\label{eqcm}
\textbf{B}(\textbf{x}) = 
\begin{cases}
b(\textbf{x}) & \text{for} \quad n = 2, \\
\big( 0,0,b(\textbf{x}) \big) & \text{for} \quad n = 3,
\end{cases}
\end{equation}  
be a magnetic field such that $b(\textbf{x}) = b(x,y)$ is 
\textit{an admissible magnetic field}. Namely, there exists 
a constant $b_{0} > 0$ satisfying
\begin{equation}
b(x,y) = b_{0} + \widetilde{b}(x,y),
\end{equation}
where $\widetilde{b}$ is such that the Poisson equation 
$\Delta \widetilde{\varphi} = \widetilde{b}$
admits a solution $\widetilde{\varphi} \in C^{2}(\mathbb{R}^{2})$ satisfying
$\sup_{(x,y) \in \mathbb{R}^{2}} \vert D^{\alpha} \widetilde{\varphi}(x,y) 
\vert < \infty$, $\alpha \in \bz_{+}^{2}$, $\vert \alpha \vert \leq 2$. 
By defining $\varphi_{0} (x,y) := \frac{1}{4}b_{0} (x^2 + y^2)$ 
on $\br^2$ and
\begin{equation}  
\varphi (x,y) := \varphi_{0} (x,y) + \widetilde{\varphi} (x,y),
\end{equation}
we obtain a magnetic potential 
$\textbf{A}_n : \br^n \longrightarrow \br^n$ generating the magnetic field 
$\textbf{B} = \text{curl} \, \textbf{A}_n $ by setting 
\begin{equation}
\textbf{A}_n(\textbf{x})= 
\begin{cases}
\textbf{A}_n(x,y) = 
\big( - \partial_y \varphi (x,y),\partial_x \varphi (x,y) \big) & \text{for} 
\quad n = 2, \\
\textbf{A}_n(x,y,\xpa) = 
\big( - \partial_y \varphi (x,y), \partial_x \varphi (x,y),0 \big) & \text{for} 
\quad n = 3.
\end{cases}
\end{equation} 

\begin{rem}
\end{rem}
\begin{itemize}
\item[(i)] The class of admissible magnetic fields described above is essentially 
the one introduced in \cite{ra,rage}. We refer to these papers for more details and 
examples of admissible magnetic fields.
\item[(ii)] In \eqref{eqcm}, $b$ stands for the intensity of the magnetic field.
\item[(iii)] The case $\widetilde{b} = 0$ corresponds to the constant magnetic 
field of strength $b_0 > 0$. 
\item[(iv)] In the three-dimensional case $n = 3$, the magnetic field is of 
constant direction and points in the $\xpa$-direction.
\end{itemize}
Let $V(\textbf{x}) = \big\lbrace V_{\ell k}(\textbf{x}) \big\rbrace_{\ell,k = 1}^2$ 
be a $2 \times 2$ complex matrix-valued potential. Then, the Pauli operators $\hn(b,V)$ 
acting on $L^{2}(\br^n) := L^{2}(\br^n,\bc^2)$, $n = 2$, $3$, are defined by 
\begin{equation}\label{eq1,1}
\hn(b,V) := \begin{pmatrix}
   (-i\nabla - \textbf{A}_n)^{2} - b & 0 \\
   0 & (-i\nabla - \textbf{A}_n)^{2} + b
\end{pmatrix} + V,
\end{equation}
initially on $C_{0}^{\infty}(\br^n,\mathbb{C}^{2})$, and then closed in 
$L^{2}(\br^n)$. 

\medskip

For $V = 0$, we have the following result from \cite[Propositions 1.1 and 1.2]{ra}  
about the spectrum $\sigma \big( \hd(b,0) \big)$ of the operator $\hd(b,0)$:

\begin{prop}\label{pr}
Let $b$ be an admissible magnetic field with $b_0 > 0$. Then, 
$0 = \inf \sigma \big( \hd(b,0) \big)$ is an isolated eigenvalue of infinite 
multiplicity. More precisely, we have
\begin{equation}
\dim \, {\rm Ker} \big( (-i\nabla - \textbf{A}_2)^{2} - b \big) = \infty, \quad
\dim \, {\rm Ker} \big( (-i\nabla - \textbf{A}_2)^{2} + b \big) = 0,
\end{equation}
and 
\begin{equation}
(0,\zeta) \subset \br \setminus \sigma \big( \hd(b,0) \big),
\end{equation}
where 
\begin{equation}\label{eq1,6}
\zeta := 2 b_{0} e^{-2 \hspace{0.5mm} \textup{osc} \hspace{0.5mm} \widetilde{\varphi}},
\qquad \textup{osc} \hspace{0.5mm} \tilde{\varphi}:= \sup_{(x,y) \in \br^{2}} 
\widetilde{\varphi} (x,y) - \inf_{(x,y) \in \br^{2}} \widetilde{\varphi} (x,y).
\end{equation}
\end{prop}

\noindent
In particular, by \cite[Corollary 2.2]{rage}, we have
\begin{equation}\label{eqsp}
\sigma \big( \htr(b,0) \big) = \sigma_{\textbf{ac}} \big( \htr(b,0) \big) = 
[0,\infty).
\end{equation} 

\noindent
Throughout this paper, our minimal assumptions on the 
potentials $V$ are the following:

\medskip

\noindent
\textbf{Assumption (A1):} 
For $n = 2$, we assume that
\begin{equation}\label{eq1,12}
\begin{aligned}
0 \not\equiv V_{\ell k}, \quad \vert V_{\ell k}(x,y) 
\vert \leq F (x,y), \quad 1 \leq \ell,k 
\leq 2,
\end{aligned}
\end{equation}
where $F \in \bigl( L^\frac{q}{2} \cap L^\infty \bigr) \big( \br^{2},
\br_{+}^{\ast} \big)$ for some $2 \leq q < \infty$. 

\medskip

\noindent
\textbf{Assumption (C1):}
For $n = 3$, we assume that
\begin{equation}\label{eq1,13}
\begin{split}
& \ast \hspace{0.6mm} 0 \not\equiv V_{\ell k}, \quad 
\vert V_{\ell k}(x,y,X_\parallel) \vert \leq G_{\perp} 
(x,y) \hspace{0.5mm} G(X_\Vert), \quad 1 \leq \ell,k \leq 2, \\
& \ast \hspace{0.6mm} G_{\perp} \in \bigl( L^\frac{q}{2} \cap L^\infty 
\bigr) \big( \br^{2},\br_{+}^{\ast} \big) \; 
\textup{for some $2 \leq q < \infty$}, \\
& \ast 0 < G(X_\Vert) \leq \text{Const.} 
\hspace{0.5mm} \langle X_\Vert \rangle^{-m}, m > 3, \hspace{0.4mm} 
\textup{where} \hspace{0.4mm} \langle y \rangle := 
(1 + \vert y \vert^2)^\frac{1}{2} \: \textup{for} \: y \in \br^d.
\end{split}
\end{equation}


\noindent
\textbf{Examples:}
\begin{itemize}
\item[(i)] In Assumptions (A1) and (C1), both $F$ and $G_{\perp}$ can be 
thought of the function
$\br^2 \ni (x,y) \mapsto \langle (x,y) \rangle^{-m_{\perp}}$ with
$m_{\perp} > 0$.
\item[(ii)] In Assumption (C1), nonreal-valued potentials $V$ with
$V_{\ell k}(x,y,X_\parallel) = \mathcal{O} \big( \langle (x,y,X_\parallel) 
\rangle^{-\alpha} \big)$, $1 \leq \ell,k \leq 2$, $\alpha > 3$,
can be considered since obviously this implies that 
$$
V_{\ell k}(x,y,X_\parallel) = \mathcal{O} \big( \langle (x,y) \rangle^{-m_{\perp}} 
\langle X_\Vert \rangle^{-m} \big), \quad m \in (3,\alpha), \quad 
m_{\perp} = \alpha - m > 0.
$$ 
\end{itemize}

\subsection{State of the article}

Since we will deal with non-self-adjoint operators, for 
convenience, we introduce some conventional definitions and 
notations. Let $S$ be a closed operator acting on a separable 
Hilbert space. An isolated point $\mu$ of $\sigma(S)$ lies in 
$\sigma_{\textup{\textbf{disc}}}(S)$, the discrete spectrum 
of $S$, if it's algebraic multiplicity
\begin{equation}\label{eq1,15}
\textup{mult}(\mu) := \textup{rank} \left( 
\frac{1}{2i\pi} \int_{\mathscr{C}} (S - z)^{-1}dz 
\right)
\end{equation}
is finite, $\mathscr{C}$ being a small positively oriented 
circle centred at $\mu$ and containing $\mu$ as the only point 
of $\sigma(S)$. Note that the geometric multiplicity of $\mu$, 
defined by $\dim \, {\rm Ker} \, (S - \mu)$, satisfies the 
inequality $\dim \, {\rm Ker} \, (S - \mu) \leq 
\text{mult}(\mu)$, equality happening if $S$ is self-adjoint. 
We define the essential spectrum 
$\sigma_{\text{\textbf{\textup{ess}}}}(S)$ of $S$ as the set 
of points $\mu \in \bc$ such that $S - \mu$ is not a Fredholm 
operator.
Under Assumptions (A1) and (C1), we prove that $V$ is relatively 
compact with respect to $\hn(b,0)$, $n = 2$, $3$. Therefore, due 
to the Weyl criterion on the invariance of the essential spectrum,
we have 
$
\sigma_{\text{\textbf{ess}}} \big( \hn(b,V) \big) = 
\sigma_{\text{\textbf{ess}}} \big( \hn(b,0) \big)
$, $n = 2$, $3$.
However, the potential $V$ may generate (complex) discrete spectrum 
whose only accumulation points are 
$\sigma_{\text{\textbf{ess}}} \big( \hn(b,V) \big)$, see 
\cite[Theorem 2.1, p. 373]{goh}. The distribution of the discrete 
spectrum near the essential spectrum for the quantum Hamiltonians
has been extensively studied by various authors. However, most of 
the known results treat the case of self-adjoint electric potentials, 
see for instance \cite[Chap. 11-12]{iv}, \cite{mel,sob,tam,roz,ra,bon,rage,bo} 
and the references given there. But, recently
and during the past years, there has been an increasing interest in 
the spectral theory of non-self-adjoint differential operators, in 
particular for the quantum Hamiltonians, see for instance \cite{fra, 
bru1,bor,dem,lap,demu,gol,han,wan,dio,dub,cue,fra1}. For a detailed 
bibliography on the theory, we refer for instance to \cite{wan,cue}. 
Another results on spectral properties for non-self-adjoint operators 
can be found in Sjöstrand paper \cite{sjo} and the references given 
there. Results concerning non-self-adjoint Pauli operators are much 
more sparse, see for instance \cite{koc}, where the authors investigated 
the 1D Pauli equation with complex boundary conditions. 

\medskip

The aim of the present paper is to describe simple methods of obtaining 
complex eigenvalues asymptotics near the low ground energy $0$ of the 2D 
Pauli operator $\hd(b,V)$, and to show how we can construct complex 
matrix-valued potentials $V$ generating nonreal eigenvalues near the low 
ground energy $0$ of the 3D Pauli operator $\htr(b,V)$. Our work is closely 
related to \cite{dio2,dio1}, where the author treats the case of the 
Schrödinger and Dirac operators with constant and nonconstant magnetic 
fields. Both in these papers and in the present one, the proofs of the 
results are inspired by previous works (on characteristic values
and resonances) for self-adjoint perturbations (see \cite{bon,bo}). 
More precisely, here, in the 2D case, the spectral gap $(0,\zeta)$ in 
$\sigma \big( \hd(b,0) \big)$ allows to reduce the study of 
$\sigma_{\text{\textbf{disc}}} \big( \hd(b,V) \big)$ near $0$, to that of the
zeros of a holomorphic function in \emph{a punctured neighbourhood of} $0$ 
(see Lemma \ref{l4,1} and Proposition \ref{p4,1}). Hence, by this way, we can 
apply the general approach developed in \cite{bo} to solve our problem. On the
contrary, in the 3D case, since $\sigma \big( \htr(b,0) \big)$ is absolutely continuous, 
this reduction holds in \emph{half-disks not containing} $0$ (see Lemma \ref{l7,3} 
and Proposition \ref{p7,2}). In that case, \cite{bo}'s approach does not work 
and we have to use the one developed in \cite{bon} to solve our problem. 
The methods of this article also combine functional analysis, complex analysis, 
functional determinant and spectral properties of Toeplitz operators which appear 
when making spectral reduction near the low ground energy $0$. The main difficulties 
come from the matrix-valuedness and the non-selfajointness of $V$. Moreover, in the 
three-dimensional case, unlike the study of the resonances, $V$ is not 
supposed to be exponentially decreasing with respect to the direction of 
the magnetic field. Thus, the operator-valued function $z \longmapsto 
\mathcal{T}_V(z)$ defined in Lemma \ref{l7,3} is not analytic near the real 
axis and some limiting absorption principle have to be used (see Proposition 
\ref{p7,2}). In contrast with the Laplace operator, in \cite{wan}, Wang 
investigated the case of $-\Delta + V$ in $L^2 \big( \br^n \big)$, $n \geq 2$, 
$V$ being a dissipative potential, $i.e.$ $V(x) = V_{1}(x) - iV_{2}(x)$ where 
$V_{1}$ and $V_{2}$ are two measurable functions satisfying $V_{2}(x) \geq 0$, 
and $V_{2}(x) > 0$ on an open non empty set. He proved that $0$ is not 
an accumulation point of the complex eigenvalues if the potential decays 
more rapidly than $|x|^{-2}$. It is still unknown, for more general 
complex potentials without sign restriction on the imaginary part, whether 
$0$ can be an accumulation point of complex eigenvalues or not. 
In the present work, we show that in the presence of a magnetic field, 
the situation is totally different. Even a compactly supported perturbation 
can produce clusters of eigenvalues near the low ground energy $0$ of the
operators $\hn(b,V)$, $n = 2$, $3$. More precisely, in the case $n = 3$, 
for some sufficiently small and sufficiently decreasing potential of the form
\begin{itemize}
\item $\eta W$, with $W$ a positive Hermitian matrix,
\item $\eta \in \bc^\ast$, $Arg(\eta) \in \pm (\frac{\pi}{2},\pi)$,
\end{itemize}
we prove that $0$ is an accumulation point of a sequence of eigenvalues 
which are concentrated along the semi-axis 
$e^{i(2Arg(\eta) \mp \pi)} [0,+\infty)$ (see Theorem \ref{t2,6}). On 
the contrary, when 
\begin{itemize}
\item $Arg(\eta) \in \pm (0,\frac{\pi}{2})$,
\end{itemize} 
this phenomenon disappears (see Corollary \ref{c1}). In the case $n = 2$,
the situation is rather different in the sense we prove that $0$ is an 
accumulation point of a sequence of eigenvalues which are concentrated 
along the semi-axis $\pm e^{iArg(\eta)} [0,+\infty)$, for some sufficiently 
decreasing potential of the form
\begin{itemize}
\item $\eta W$, with $\pm W$ a positive Hermitian matrix,
\item $\eta \in \bc^\ast$,
\end{itemize}
(see Theorem \ref{t2,2}). The case of the Laplace 
operator is also studied in a recent preprint by Bögli \cite{bog}, where nonreal 
potentials decaying at infinity generating infinitely many nonreal eigenvalues 
accumulating at each point of the essential spectrum $[0,+\infty)$ are 
constructed. Upper bounds on the number of the complex eigenvalues in small 
annulus near the low ground energy $0$ of the operators $\hn(b,V)$, $n = 2$, 
$3$, are also established here (see Theorems \ref{t2,1} and \ref{t2,4} respectively).

\subsection{Organisation of the paper}

The paper is organized as follows. Our main results are stated in Section 
\ref{s2}. Section \ref{s4} is devoted to the study of the discrete spectrum 
near the low ground energy for the two-dimensional Pauli operator. The 
corresponding main results are proved in Sections \ref{s5} and \ref{s6}. 
Section \ref{s7} is devoted to the study of the discrete spectrum near the 
low ground energy with respect to the three-dimensional Pauli operator, 
the corresponding main results being proved in Sections \ref{s8} and \ref{s10}. 
Section \ref{s3,1} is a brief appendix on basic properties of Schatten-von 
Neumann class ideals, Section \ref{sa,1} a brief appendix on the theory of 
the index of a finite meromorphic operator-valued function, and Section 
\ref{sa,3} a brief appendix on the notion of characteristic values of 
operator-valued functions.

\subsection*{Notations}

For a $2 \times 2$ matrix $M : \bc^2 \rightarrow \bc^2$, 
$\vert M \vert$ denotes the multiplication operator by the 
matrix
$\sqrt{M^\ast M}(\textbf{x}) =: \big\lbrace \vert M \vert_{\ell k}
(\textbf{x}) \big\rbrace$, $1 \leq \ell,k \leq 2$, $\textbf{x} 
\in \mathbb{R}^{n}$, $n = 2, \hspace{0.6mm} 3$.
We will denote $\mathfrak{B}_h(\br^n)$ the set of $2 \times 2$ Hermitian 
matrices on $\br^n$, $n = 2$, $3$. The spectral projection of $L^2 (\br^2)$ 
onto the (infinite-dimensional) kernel of
$\hd^- := (-i\partial_x - a_1)^2 + (-i\partial_y - a_2)^2 - b$,
will be denoted $p := p(b)$. Here, $a_j$, $j = 1$, $2$, are the components of
the magnetic potential ${\bf A}_2$, so that $\hd^-$ is the first component of
the operator $\hd(b,0)$. The operator 
$\hd^+ := (-i\partial_x - a_1)^2 + (-i\partial_y - a_2)^2 + b$ will denote 
his second component.

\section{Statement of the main results}\label{s2}

This section is devoted to the formulation of our main results. 
The eigenvalues will be counted according to their algebraic 
multiplicity defined above. As preparation, we first recall some 
well-known results on Toeplitz operators. We know from \cite[Lemma 2.3]{rage} that 
if $U \in L^q(\br^2)$, $q \geq 1$, then the Toeplitz operator $pUp$ 
belongs to the Schatten-von Neumann class $\sqq \big( L^2(\br^2) \big)$ 
(see Section \ref{s3,1} for the definition of the Schatten classes 
$\sqq$). In particular, $pUp$ is a compact operator. Moreover, when 
it is self-adjoint and positive, the following asymptotics about the 
quantity $\textup{Tr} \, \textbf{\textup{1}}_{(r,\infty)} (pUp)$, 
$r \searrow 0$, are well-known:

\begin{itemize}
\item[\bf (H1)] If $0 \leq U \in C^{1} \big( \mathbb{R}^{2} \big)$ 
verifies
$U(x,y) = u_{0} \left( \frac{(x,y)}{\Vert (x,y) \Vert} \right) \Vert 
(x,y) \Vert^{-m} ( 1 + o(1) \big)$, $\Vert (x,y) \Vert \rightarrow \infty$, 
$m > 0$ constant, where $u_{0}$ is a non-negative continuous function on 
$\mathbb{S}^{1}$ not vanishing identically,
$\vert \nabla U(x,y) \vert \leq C_{1} \langle (x,y) \rangle^{-m-1}$
with some constant $C_{1} > 0$, and if there exists
an integrated density of states for the operator $P_2^-$ (see \cite{rage}
definition (3.11)), then by \cite[Lemma 3.3]{rage},
\begin{equation}\label{eq2,02}
\textup{Tr} \, \textbf{\textup{1}}_{(r,\infty)} 
\big( p U p \big) = 
C_{m} r^{-2/m} \big( 1 + o(1) \big), \hspace{0.2cm} r \searrow 0,
\end{equation}
where
$C_{m} := \frac{b_0}{4\pi} \int_{\mathbb{S}^{1}} dt \, u_{0}(t)^{2/m}$.

\item[\bf (H2)] If $0 \leq U \in L^{\infty}\big(\mathbb{R}^{2}\big)$ 
verifies 
$\ln U(x,y) = -\mu \Vert (x,y) \Vert^{2\beta} \big( 1 + o(1) \big)$, 
$\Vert (x,y) \Vert \rightarrow \infty$,
with some constants $\beta > 0$ and $\mu > 0$, then by 
\cite[Lemma 3.4]{rage},
\begin{equation}\label{eq2,3}
\textup{Tr} \, \textbf{\textup{1}}_{(r,\infty)} \big( p U p \big) = 
\varphi_{\beta}(r) \big( 1 + o(1) \big), \hspace{0.2cm} r \searrow 0,
\end{equation}
where for $0 < r < \textup{e}^{-1}$, we set
$$
\varphi_{\beta}(r) :=
 \begin{cases}
 \frac{1}{2} b_0 \mu^{-1/\beta} \vert \ln r \vert^{1/\beta} & \text{if } 
 0 < \beta < 1,\\
 \frac{1}{\ln(1 + 2\mu/b_0)} \vert \ln r \vert & \text{if } \beta = 1,\\
 \frac{\beta}{\beta - 1} \big( \ln \vert \ln r \vert \big)^{-1} \vert 
 \ln r \vert & \text{if } 1 < \beta < \infty.
 \end{cases}
$$

\item[\bf (H3)] If $0 \leq U \in L^{\infty} \big( \mathbb{R}^{2} \big)$ 
has a compact support and if there exists a constant $C > 0$ verifying
$C \leq U$ on an open non-empty subset of $\mathbb{R}^{2}$, then by 
\cite[Lemma 3.5]{rage},
\begin{equation}\label{eq2,4}
\textup{Tr} \hspace{0.4mm} \textbf{\textup{1}}_{(r,\infty)} 
\big( p U p \big) = 
\varphi_{\infty}(r) \big( 1 + o(1) \big), 
\hspace{0.2cm} r \searrow 0,
\end{equation}
where 
$
\varphi_{\infty}(r) := \big( \ln \vert \ln r \vert \big{)}^{-1} 
\vert \ln r \vert, \hspace{0.2cm} 0 < r < \textup{e}^{-1}
$.
\end{itemize}

\subsection{Results for the case 2D}\label{s2,1}

Let $V$ satisfy Assumption (A1). In that case, we have $\vert V \vert_{\ell k}(x,y) 
= \mathcal{O} \big( F(x,y) \big)$ for any $1 \leq \ell,k \leq 2$ and any 
$(x,y) \in \br^2$. Then, \cite[Lemma 2.3]{rage} implies that the Toeplitz operator 
$p \vert V \vert_{11} p$ belongs to the Schatten-von Neumann class 
$\mathcal{S}_{\frac{q}{2}} \big( L^2(\br^2) \big)$, $q \geq 2$. In particular,
$p \vert V \vert_{11} p$ is a self-adjoint positive compact operator. For $\zeta$ 
defined by \eqref{eq1,6}, fix $0 < \epsilon \leq \zeta$ and introduce the punctured 
disk of $0$
\begin{equation}\label{eq2,5}
D(0,\epsilon)^{\ast} := \big\lbrace \mu \in \mathbb{C} : 0 < \vert \mu \vert < 
\epsilon \big\rbrace.
\end{equation}
Our first main result gives an upper bound on the number of 
complex eigenvalues near the low ground energy $0$ of the 
operator $\hd(b,0)$, in small annulus.

\begin{theo}[Local upper bound]\label{t2,1}
Assume that Assumption (A1) holds with $\zeta > \Vert V \Vert \ll 1$. 
Then, there exists $0 < r_0 < \zeta$ such that for any $r > 0$ with  
$r < r_{0} < \frac{3}{2}r$,
\begin{equation}\label{eq2,6}
\# \Big\lbrace \mu \in \sigma_{\textup{\textbf{disc}}} \big( \hd(b,V) \big) : 
r < \vert \mu \vert < 2r \Big\rbrace \leq C \, \textup{Tr} \,
\one_{(r,\infty)} \big( p \vert V \vert_{11} p \big) 
\vert \ln r \vert + \mathcal{O}(1),
\end{equation}
for some constant $C > 0$. In particular, if $\vert V \vert_{11}$ is compactly
supported, then $\textup{Tr} \, \one_{(r,\infty)} \big( p \vert V \vert_{11} p 
\big)$ satisfies the asymptotic \eqref{eq2,4} as $r \searrow 0$.
\end{theo}

\noindent
In order to get asymptotic results, we put some restrictions on 
the potential $V$.

\medskip

\noindent
{\bf Assumption (A2):}
We assume that
\begin{equation}\label{eq2,7}
\small{V = \eta W \hspace{0.1cm} \text{with} \hspace{0.1cm} \eta \in \bc^\ast,
\hspace{0.1cm} \text{and} \hspace{0.1cm}
W(x,y) = \big\lbrace W_{\ell k}(x,y) \big\rbrace_{\ell,k = 1}^2 \in 
\mathfrak{B}_h(\br^2), \hspace{0.1cm} (x,y) \in \mathbb{R}^{2}.}
\end{equation} 

\noindent
Denote $J := sign(W)$ the matrix sign of $W$ satisfying 
$W = J \vert W \vert$. Let $\text{P}_\perp$ and 
$\text{Q}_\perp$ be the orthogonal projections in 
$L^{2} \big( \br^{2} \big)$ defined by
\begin{equation}\label{eq2,8}
\textup{P}_\perp := \begin{pmatrix} 
   p & 0 \\ 
   0 & 0 
\end{pmatrix}, \hspace{1cm} \textup{Q}_\perp := 
\textup{I} - \textup{P}_\perp = \begin{pmatrix} 
   I - p & 0 \\
   0 & I 
\end{pmatrix}.
\end{equation}
For $\mu \in D (0,\epsilon) := D (0,\epsilon)^{\ast} 
\cup \lbrace 0 \rbrace$, we introduce the operator 
\begin{equation}\label{eq2,9}
\small{\displaystyle \mathcal{A}_\perp(\mu) := 
J \vert W \vert^{\frac{1}{2}} p \begin{pmatrix} 
   1 & 0 \\ 
   0 & 0 
\end{pmatrix} \vert W \vert^{\frac{1}{2}} -
\mu J \vert W \vert^{\frac{1}{2}} 
\big( \hd(b,0) - \mu \big)^{-1} \textup{Q}_\perp
 \vert W \vert^{\frac{1}{2}}.}
\end{equation}
Note that $\mathcal{A}_\perp(\cdot)$ is 
holomorphic on $D (0,\epsilon)$. 
Introduce $\Pi$ the orthogonal projection onto 
$\text{Ker} \, \mathcal{A}_\perp(0)$ together with 
the following condition:
\begin{equation}\label{eq2,10}
\begin{cases}
\textup{$I - \eta \mathcal{A}_\perp '(0)\Pi$ is an invertible operator for} \\
\textup{$W$ of definite sign $\pm W(x,y) \geq 0$, $(x,y) \in \br^2$.}
\end{cases}
\end{equation}

\begin{rem}\label{r2,1}
\end{rem}
\begin{itemize}
\item[(i)] There is no loss of generality in saying $W$ is of 
definite sign $J = \pm$.
\item[(ii)] Actually, $\mathcal{A}_\perp '(0)\Pi$ is compact so 
that condition \eqref{eq2,10} is generically satisfied. Namely, it 
is fulfilled whenever $\eta \in \br \setminus \lbrace e_n^{-1}, n \in 
\bn \rbrace$, where $(e_n)_n$ denotes the nonzero eigenvalue sequence 
of the operator $\mathcal{A}_\perp '(0)\Pi$.
\item[(iii)] Condition \eqref{eq2,10} is also fulfilled for 
small potentials $V$. Namely, for potentials $V = \varepsilon 
\eta W$ with $\varepsilon$ a real number small enough.
\end{itemize}

\noindent
If $r_{0}$, $\delta$, are two positive fixed constants, 
and $r > 0$ tending to zero, we define the sector
\begin{equation}\label{eq2,11}
\Gamma^{\delta}(r,r_{0}) := \big\lbrace x + iy \in \mathbb{C} 
: r < x < r_{0}, -\delta x < y < \delta x \big\rbrace.
\end{equation}

\begin{theo}[Localization, asymptotic behaviours]\label{t2,2}
Let $V$ satisfy Assumptions (A1) and (A2), with $\eta \in \bc^\ast$
and $W$ of definite sign $J = \pm$. Assume that \eqref{eq2,10}
holds. Then, there exits $r_{0} > 0$ such that near zero:

\begin{itemize}
\item[(i)] Localization: The discrete eigenvalues $\mu$ of 
$\hd(b,\eta W)$ with $0 < \vert \mu \vert < r_{0}$ satisfy
\begin{equation}\label{eq2,12}
\mu \in \pm \eta  \hspace{0.1cm} 
\overline{\Gamma^{\delta}(r,r_{0})},
\end{equation}
for any $\delta > 0$.
\item[(ii)] Asymptotic: There exists a sequence 
$(r_{\ell})_{\ell}$ of positive numbers tending to zero such 
that
\begin{equation}\label{eq2,13}
\# \Big\lbrace \sigma_{\textup{\textbf{disc}}} \big( \hd(b,\eta W) \big) : 
r_{\ell} < \vert \mu \vert < r_0 \Big\rbrace = \textup{Tr} \hspace{0.4mm} 
\one_{(r_{\ell},\infty)} \big( p \vert W \vert_{11} p \big) \big( 1 + o(1) \big),
\end{equation}
as $\ell \longrightarrow \infty$.
\item[(iii)] Asymptotic: If $\vert W \vert_{11}$ satisfies the hypotheses
(H1), (H2) and (H3) above, then,
\begin{equation}\label{eq2,130}
\# \Big\lbrace \sigma_{\textup{\textbf{disc}}} \big( \hd(b,\eta W) \big) : 
r < \vert \mu \vert < r_0 \Big\rbrace = \textup{Tr} \,
\one_{(r,\infty)} \big( p \vert W \vert_{11} p \big) \big( 1 + o(1) \big),
\end{equation}
as $r \searrow 0$.
\end{itemize}
\end{theo}


\begin{figure}[h]\label{fig 2}
\begin{center}
\tikzstyle{+grisEncadre}=[fill=gray!60]
\tikzstyle{blancEncadre}=[fill=white!100]
\tikzstyle{grisEncadre}=[fill=gray!25]
\tikzstyle{dEncadre}=[dotted]

\begin{tikzpicture}[scale=0.8]

\draw (0,0) circle(2);


\draw [+grisEncadre] (0,0) -- (39:2) arc (39:73:2) -- cycle;

\draw (0,0) circle(0.5);

\draw[->] [thick] (-2.6,0) -- (3,0);
\draw (3,0) node[right] {\tiny{$\Re(\mu)$}};

\draw[->] [thick] (0,-2.5) -- (0,2.5);
\draw (0,2.5) node[above] {\tiny{$\Im(\mu)$}};

\draw (-0.15,0.02) node[above] {\tiny{$r$}};
\draw (0,0) -- (-0.45,0.17);

\draw (0,0) -- (1.7,2.5);
\draw (2.4,2.3) node[above] {$y = \tan \big( Arg(\eta) \big) x$};

\draw (0,0) -- (0.8,-1.84);
\draw (0.6,-1.2) node[above] {\tiny{$r_{0}$}};

\draw (0.95,1.7) node[above] {\tiny{$\theta$}};
\draw (1.4,1.4) node[above] {\tiny{$\theta$}};

\node at (0.45,0.5) {\tiny{$\times$}};
\node at (0.45,0.65) {\tiny{$\times$}};
\node at (0.63,0.9) {\tiny{$\times$}};
\node at (0.76,1.1) {\tiny{$\times$}};
\node at (0.8,0.95) {\tiny{$\times$}};
\node at (0.62,0.75) {\tiny{$\times$}};

\node at (0.3,0.6) {\tiny{$\times$}};
\node at (0.4,0.8) {\tiny{$\times$}};
\node at (0.6,1.1) {\tiny{$\times$}};
\node at (0.5,0.95) {\tiny{$\times$}};

\node at (0.35,0.35) {\tiny{$\times$}};
\node at (0.2,0.3) {\tiny{$\times$}};
\node at (0.25,0.45) {\tiny{$\times$}};
\node at (0.6,0.6) {\tiny{$\times$}};
\node at (0.8,0.8) {\tiny{$\times$}};
\node at (1.2,1.4) {\tiny{$\times$}};
\node at (1,1.6) {\tiny{$\times$}};
\node at (0.7,1.7) {\tiny{$\times$}};
\node at (1,1.2) {\tiny{$\times$}};
\node at (1.1,1) {\tiny{$\times$}};

\node at (-4.2,2.5) {$V = \eta W$};
\node at (-4.2,1.8) {$\eta \in \bc^\ast, \hspace{0.5mm} W \geq 0$};

\end{tikzpicture}
\caption{\textbf{A graphic illustration of the localization of 
the nonreal eigenvalues near zero:} \textit{For $r_{0}$ small 
enough, the nonreal 
eigenvalues $\mu$ of $\hd(b,\eta W)$ are localized near
the semi-axis $\mu = \eta ]0,+\infty)$ in small angular 
sectors.}}
\end{center}
\end{figure}
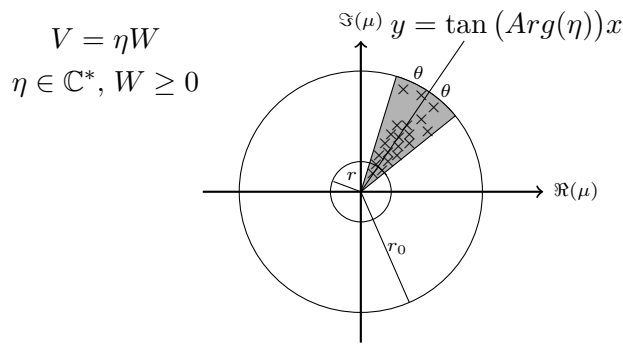

\noindent
\textbf{Corollaries and Remarks.}
\begin{itemize}
\item[(i)] According to \eqref{eq2,02}, \eqref{eq2,3} and \eqref{eq2,4}, 
the asymptotics obtained in Theorem \ref{t2,2} essentially coincide
with those obtained in \cite{ra}, where the case of the self-adjoint Pauli 
operator $\hd(b,V)$ with more general assumptions on $V$ is considered. Novelty 
in this paper is that we consider complex matrix-valued potentials $V$.

\item[(ii)] Theorem \ref{t2,2} is still true if we 
replace the assumption $J = \pm$ by 
$$
sign(W) \begin{pmatrix} 
   1 & 0 \\ 
   0 & 0 
\end{pmatrix} = \begin{pmatrix} 
   \pm 1 & 0 \\ 
   0 & 0 
\end{pmatrix}.
$$
This condition is for instance fulfilled by potentials 
$W \in \mathfrak{B}_h(\br^2,\bc^2)$ of the form 
$$
W = \begin{pmatrix} 
  \pm W_{11} & 0 \\ 
   0 & W_{22} 
\end{pmatrix}, \quad W_{11} > 0.
$$ 

\item[(iii)] The existence of nonreal eigenvalues and 
their accumulation near $0$ are ensured by assertion (ii) 
and (iii) of Theorem \ref{t2,2}. 

\item[(iv)] Another way of stating assertion (i) of 
Theorem \ref{t2,2} is as follows: near the low 
ground energy of $\hd(b,0)$, the discrete eigenvalues $\mu$ 
of $\hd(b,\eta W)$ are concentrated around the semi-axis 
$\pm \eta ]0,+\infty)$, respectively for $J = \pm$. 

\item[(v)] For $\eta = 1$, we recover the standard 
self-adjoint electric potentials $V$ of definite sign, and 
\eqref{eq2,12} becomes $\mu \in \pm \hspace{0.1cm} ]0,+\infty)$. 
This just amounts to saying that the discrete eigenvalues 
are localized on the right and the left of $0$. 

\item[(vi)] If $\eta = i$, we can write \eqref{eq2,12} 
as $\pm \Im(\mu) \geq 0$ with $\vert \Re(\mu) \vert = o(\vert 
\mu \vert)$, see \eqref{eq6,1}. This means that the discrete 
eigenvalues are concentrated in a vicinity of the semi-axis 
$\pm i ]0,+\infty)$.
\end{itemize}

\subsection{Results for the case 3D}\label{s2,2}

Let $V$ satisfy Assumption (C1) and $\textbf{\textup{V}}_{11}$ 
be the multiplication operator by the function (also noted) 
$\textbf{\textup{V}}_{11} : \br^2 \longrightarrow \br$ 
defined by 
\begin{equation}\label{eq2,14}
\displaystyle \textbf{\textup{V}}_{11}(x,y) := \frac{1}{2}
\int_\br dX_\Vert \, \vert V \vert_{11} (x,y,X_\Vert).
\end{equation}
Then, as in the case 2D, the Toeplitz operator 
$p \textbf{\textup{V}}_{11} p$ satisfies $p \textbf{\textup{V}}_{11} 
p \in \sqq \big( L^2(\br^3) \big)$ so that it is a self-adjoint 
positive compact operator. In the sequel,
\begin{equation}\label{eq2,2}
\bc_\pm := \big\lbrace z \in \bc : \pm \Im(z) > 0 \big\rbrace
\end{equation}
denotes respectively the upper and lower half-plane. For 
$0 \leq \varrho_1 < \varrho_2 \leq \zeta$, we introduce the ring
\begin{equation}\label{eq2,3bb}
D(\varrho_1,\varrho_2) := \big\lbrace z \in \bc: 
\varrho_1 < \vert z \vert < \varrho_2  \big\rbrace,
\end{equation}
and the half-rings
\begin{equation}\label{eq2,3b}
D_\pm(\varrho_1,\varrho_2) := D(\varrho_1,\varrho_2) 
\cap \bc_\pm.
\end{equation}
For $\nu > 0$ constant, we define the domains
\begin{equation}\label{eq2,3bbb}
D_{\pm}^{\nu}(\varrho_1,\varrho_2) := D_\pm(\varrho_1,\varrho_2) 
\cap \big\lbrace z \in \bc : \vert \Im(z) \vert > \nu \big\rbrace.
\end{equation}
Our first main result gives upper bounds on the number of 
complex eigenvalues near the low ground energy $0$ of the 
operator $\htr(b,0)$, in small half-rings. 

\begin{theo}[Local upper bounds]\label{t2,4}

Assume that \textit{Assumption (C1)} holds with $\zeta > \Vert V \Vert \ll 1$. 
Then, there exists $0 < r_{0} < \sqrt{\zeta}$ such that for any 
$r > 0$ with $r < r_{0} < \sqrt{\frac{5}{2}}r$, and any $0 < \nu < 2r^2$,
\begin{equation}\label{eq2,19}
\# \Big\lbrace z \in \sigma_{\textup{\textbf{disc}}} 
\big( \htr(b,V) \big) \cap D_{\pm}^\nu(r^2,4r^2)
\Big\rbrace \leq C \, \Big( \textup{Tr} \,
\one_{(r,\infty)} \big( p \textbf{\textup{V}}_{11} 
p \big) \vert \ln r \vert \Big) + \mathcal{O}(1),
\end{equation}
for some constant $C > 0$. In particular, if 
$\textbf{\textup{V}}_{11}$ is compactly supported, then 
$\textup{Tr} \, \one_{(r,\infty)} \big( p \textbf{\textup{V}}_{11} 
p \big)$ satisfies the asymptotic \eqref{eq2,4} as $r \searrow 0$.
\end{theo}

\noindent
In what follows below, two kinds of assumption on $V$ 
are needed in additional.

\medskip

\noindent
{\bf Assumption (C2):}
We assume that
\begin{equation}\label{eq2,20}
\small{V = \eta W \hspace{0.1cm} \text{with} \hspace{0.1cm} 
\eta \in \bc \setminus \br, \hspace{0.1cm} \text{and} 
\hspace{0.1cm} W(X) = \big\lbrace W_{\ell k}(X) 
\big\rbrace_{\ell,k = 1}^2 \in \mathfrak{B}_h(\br^3), 
\hspace{0.1cm} X \in \mathbb{R}^{3}.}
\end{equation}

\begin{rem}\label{r2,2}
\end{rem}
\begin{itemize}
\item[(i)] In \eqref{eq2,20}, when $W$ is of definite sign 
($i.e.$ $\pm W \geq 0$), since the change of the sign 
consists to replace $Arg(\eta)$ by $Arg(\eta) + \pi$, then 
it is sufficient to consider only $W \geq 0$.
\item[(ii)] For $\pm \sin \big( Arg(\eta) \big) > 0$ and $W \geq 0$, the 
discrete eigenvalues $z$ of the operator $\htr(b,\eta W)$ satisfy 
$\pm \Im(z) \geq 0$. 
\end{itemize}


\noindent
With respect to Assumption (C2) above, we introduce
$\textbf{\textup{W}}_{11}$, the multiplication operator 
by the function (also noted) $\textbf{\textup{W}}_{11} : 
\br^2 \longrightarrow \br$, defined as in \eqref{eq2,14} 
with $\vert V \vert_{11}$ replaced by $\vert W \vert_{11}$.
Hence, we introduce the following exponential decay 
assumption:

\medskip

\noindent
{\bf Assumption (C3):}
We assume that the function $\textbf{\textup{W}}_{11}$ satisfies
\begin{equation}\label{eq2,25}
0 < \textbf{\textup{W}}_{11}(x,y) \leq 
e^{-C \langle (x,y) \rangle^2}, \quad C > 0.
\end{equation}

\noindent
For $\alpha \in \br$ and $\theta > 0$, we introduce the sector
\begin{equation}
E_{\pm}(\alpha,\theta) := D(0,\zeta) 
\setminus \left( e^{i (2\alpha \mp \pi)} e^{i (-2\theta,2\theta)} 
(0,\zeta) \right),
\end{equation}
where we have just excluded around the semi-axis 
$z = e^{i(2\alpha \mp \pi)} (0,\zeta)$, an angular sector of 
amplitude $4\theta$. Hence, we can formulate our second main 
result in this part as follows:

\begin{theo}[Sector free of complex eigenvalues, lower bounds]\label{t2,6}
Suppose that $V$ satisfies Assumption (C1), and Assumption (C2) with 
$W \geq 0$, $\pm Arg(\eta) \in (0,\pi)$. Then, for any $\theta > 0$
small enough, there exists $\varepsilon_0 > 0$ such that: 

\begin{itemize}
\item[(i)] For any $0 < \varepsilon \leq \varepsilon_0$,
the operator $\htr(b,\varepsilon V)$ has no discrete eigenvalues in
the sector
\begin{equation}\label{e2,26}
D_\pm(r^2,r_0^2) \cap E_{\pm}(Arg(\eta),\theta), \quad 0 \le r < r_0 < \sqrt{\zeta}.
\end{equation}
\item[(ii)] If furthermore $\textbf{\textup{W}}_{11}$ satisfies
Assumption (C3), then for any $0 < \varepsilon \leq \varepsilon_0$, 
there is an accumulation of nonreal eigenvalues of 
$\htr(b,\varepsilon V)$ near zero, in a sector around 
the semi-axis $e^{i(2 Arg(\eta) \mp \pi)} (0,+\infty)$, 
for
\begin{equation}
Arg(\eta) \in \pm \left( \frac{\pi}{2},\pi \right).
\end{equation}
More precisely, there exists a decreasing sequence of 
positive numbers $(r_\ell)_\ell$, $r_\ell \searrow 0$,
satisfying
\begin{equation}\label{lb1}
\begin{split}
\# \Biggl\lbrace z \in & \sigma_{\textup{\textbf{disc}}} \big( 
\htr(b,\varepsilon V) \big) \cap D_{\pm}(\varepsilon^2 
r_{\ell + 1}^2,\varepsilon^2 r_{\ell}^2) \cap 
\left( e^{i (2 Arg(\eta) \mp \pi)} e^{i (-2\theta,2\theta)} 
(0,\zeta) \right) \Biggr\rbrace \\
& \geq \textup{Tr} \hspace{0.4mm} \one_{(r_{\ell +1},r_\ell)} 
\big( p \textbf{\textup{W}}_{11} p \big).
\end{split}
\end{equation}
\end{itemize}
\end{theo}

\bigskip

\noindent
A graphic illustration of Theorem \ref{t2,6} is given in Figure 2.2.

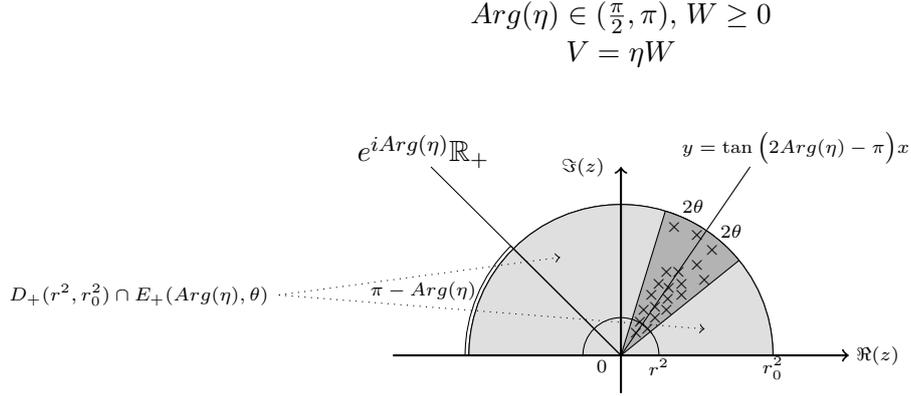
\begin{figure}\label{fig 2}
\begin{center}
\tikzstyle{+grisEncadre}=[fill=gray!60]
\tikzstyle{blancEncadre}=[fill=white!100]
\tikzstyle{grisEncadre}=[fill=gray!25]
\tikzstyle{dEncadre}=[dotted]

\begin{tikzpicture}[scale=1]

\draw (0,0) -- (0:2) arc (0:180:2) -- cycle;

\node at (-0.25,-0.15) {\tiny{$0$}};
\node at (0.5,-0.16) {\tiny{$r^2$}};
\node at (2,-0.17) {\tiny{$r_0^2$}};

\draw [grisEncadre] (0,0) -- (0:2) arc (0:180:2) -- cycle;


\draw [+grisEncadre] (0,0) -- (39:2) arc (39:73:2) -- cycle;

\draw (0,0) -- (0:0.5) arc (0:180:0.5) -- cycle;

\draw (0,0) -- (-2.5,2.5) -- cycle;
\draw (-2.6,2.35) node[above] {$e^{i Arg(\eta)} \br_+$};
\draw (-2.6,0.6) node[above] {\tiny{$\pi - Arg(\eta)$}};
\draw (0,0) -- (135:2.05) arc (135:180:2.05) -- cycle;

\draw[->] [thick] (-3,0) -- (3,0);
\draw (2.96,0) node[right] {\tiny{$\Re(z)$}};

\draw[->] [thick] (0,-0.5) -- (0,2.5);
\node at (-0.5,2.5) {\tiny{$\Im(z)$}};

\draw (0,0) -- (1.7,2.5);
\draw (2.3,2.4) node[above] {\tiny{$y = \tan \big( 2 Arg(\eta) - \pi \big) x$}};

\draw (0.95,1.75) node[above] {\tiny{$2\theta$}};
\draw (1.45,1.42) node[above] {\tiny{$2\theta$}};

\draw [dEncadre] [->] (-4.5,0.8) -- (-0.8,1.3);
\draw [dEncadre] [->] (-4.5,0.8) -- (1.1,0.35);
\draw (-4.5,0.8) node[left] {\tiny{$D_+(r^2,r_0^2) \cap E_+(Arg(\eta),\theta)$}};

\node at (0.45,0.5) {\tiny{$\times$}};
\node at (0.45,0.65) {\tiny{$\times$}};
\node at (0.63,0.9) {\tiny{$\times$}};
\node at (0.76,1.1) {\tiny{$\times$}};
\node at (0.8,0.95) {\tiny{$\times$}};
\node at (0.62,0.75) {\tiny{$\times$}};

\node at (0.3,0.6) {\tiny{$\times$}};
\node at (0.4,0.8) {\tiny{$\times$}};
\node at (0.6,1.1) {\tiny{$\times$}};
\node at (0.5,0.95) {\tiny{$\times$}};

\node at (0.35,0.35) {\tiny{$\times$}};
\node at (0.2,0.3) {\tiny{$\times$}};
\node at (0.25,0.45) {\tiny{$\times$}};
\node at (0.6,0.6) {\tiny{$\times$}};
\node at (0.8,0.8) {\tiny{$\times$}};
\node at (1.2,1.4) {\tiny{$\times$}};
\node at (1,1.6) {\tiny{$\times$}};
\node at (0.7,1.7) {\tiny{$\times$}};
\node at (1,1.2) {\tiny{$\times$}};
\node at (1.1,1) {\tiny{$\times$}};

\node at (0,4) {$V = \eta W$};
\node at (0,4.48) {$Arg(\eta) \in (\frac{\pi}{2},\pi), \hspace{0.5mm} W \geq 0$};

\end{tikzpicture}
\caption{\textbf{A graphic illustration of the localization of the nonreal 
eigenvalues near the low ground energy $0$:} 
For $\theta$ small enough and $0 < \varepsilon \leq \varepsilon_0$, 
$\htr(b,\varepsilon V)$ has no complex eigenvalues in 
$D_+(r^2,r_0^2) \cap E_+(Arg(\eta),\theta)$ (see (i)
of Theorem \ref{t2,6}). They are localized around the semi-axis 
$z = e^{i(2 Arg(\eta) - \pi)} ]0,+\infty)$ (see (ii) of 
Theorem \ref{t2,6}).}
\end{center}
\end{figure}

\begin{rem}
\end{rem}

\begin{itemize}
\item[\text{(i)}] Theorem \ref{t2,6} and \cite[Theorem 2]{bon} are quite similar in there 
structure, although in \cite{bon}, is considered the 3D Schrödinger operator with constant 
magnetic field and self-adjoint potentials decaying exponentially in the direction of the 
magnetic field, to study the resonances near the Landau levels.
\item[\text{(ii)}]If we set $\delta := \tan(\theta)$, then the proof of Theorem \ref{t2,6} 
\big(see \eqref{eq10,40} and \eqref{eq7,220}\big) shows that for $0 < r_0 < \sqrt{\zeta}$, 
the parameter $\varepsilon_0 > 0$ above depends on $\theta$ as follows:
\begin{equation}
\varepsilon_0 < \left( C \sqrt{1 + \delta^{-2}} \right)^{-1} \min \Biggl(1, 
C_1(\delta,\nu)^{-1} e^{-\Gamma_q \big( C_1(\delta,\nu) + 1 \big)^{\lceil q \rceil}} 
\Biggr),
\end{equation}
for some uniform constants $C$, $\nu$ positive, where $C_1(\delta,\nu)$ is
the constant defined by \eqref{eq7,221}, $\Gamma_q$ by \eqref{eq3,5}, and 
$\lceil q \rceil$ by \eqref{eq3,2}.
\end{itemize}

\noindent
An immediate consequence of assertion (i) of Theorem \ref{t2,6} together 
with (ii) of Remark \ref{r2,2} is the following:

\begin{cor}[Non cluster of complex eigenvalues]\label{c1}
Let the assumptions of Theorem \ref{t2,6} be fulfilled. 
Then, for any $\eta$ satisfying
$Arg(\eta) \in \pm \left( 0,\frac{\pi}{2} \right)$,
there is no accumulation of discrete eigenvalues of 
$\htr(b,\varepsilon V)$ near zero, for 
$0 < \varepsilon \leq \varepsilon_0$.
\end{cor}

\noindent
Actually, we expect this to be a general phenomenon in the following 
sense:

\begin{cnj}\label{c2,2}
Let $V = \eta W$ satisfy Assumption (C1), with $\eta \in \bc \setminus 
\br e^{ik \frac{\pi}{2}}$, $k \in \bz$, and $W \in 
\mathfrak{B}_h(\br^3)$ of definite sign. Then, there is no 
accumulation of complex eigenvalues of $\htr(b,V)$ near zero if and only if
$\Re(V) > 0$.
\end{cnj}

\begin{rem}
\textup{To make our results in perspective, let us mention that more general 
situations such as: more general magnetic field in the 3D case (not supported 
only in the $\xpa$ axis), and/or to consider perturbations of the magnetic field 
itself, are open problems.}
\end{rem}

\section{Discrete eigenvalues for the 2D problem}\label{s4}

Throughout Sections \ref{s4}-\ref{s6}, $D(0,\epsilon)^\ast$ 
stands for the punctured disk given by \eqref{eq2,5}. Here and 
in the rest of the paper, $\zeta$ is the constant defined by 
\eqref{eq1,6}.

\subsection{Preliminary results}

Let $\text{P}_\perp$ and $\text{Q}_\perp := I - \text{P}_\perp$ be 
the orthogonal projections defined by \eqref{eq2,8}. 
For $\mu \notin \sigma \big( \hd(b,0) \big)$, on account of 
\eqref{eq1,1} with $n = 2$ and Proposition \ref{pr}, we clearly have
\begin{equation}\label{eq4,1}
\big( \hd(b,0) - \mu \big)^{-1} \text{P}_\perp = 
-p \begin{pmatrix}
\mu^{-1} & 0 \\
   0 & 0
\end{pmatrix}.
\end{equation}
Therefore, for any $\mu$ lying in the resolvent set of 
$\hd(b,0)$, we have
\begin{equation}\label{eq4,2}
\big( \hd(b,0) - \mu \big)^{-1} = 
-p \begin{pmatrix}
\mu^{-1} & 0\\
0 & 0
\end{pmatrix} + \big( \hd(b,0) - \mu \big)^{-1} \text{Q}_\perp.
\end{equation}

\noindent
We begin with a general result on the first term of the r.h.s. 
of \eqref{eq4,2}.

\begin{lem}\label{l4,1} 
Let $U \in L^q(\br^2)$ with
$q \in [2,+\infty)$. Then, the operator-valued function 
$$
D(0,\epsilon)^{\ast} \ni \mu \longmapsto U 
\big( \hd(b,0) - \mu \big)^{-1} \textup{P}_\perp 
$$ 
is holomorphic with values in $\sqq \big( L^{2}(\br^2) \big)$.
Furthermore,
\begin{equation}\label{eq4,3}
\left\Vert U 
\big( \hd(b,0) - \mu \big)^{-1} \textup{P}_\perp \right\Vert_\sqq^q \leq
\frac{b_0 e^{2\textup{osc} \hspace{0.5mm} \tilde{\varphi}}}{2\pi\mu^q}
\Vert U \Vert_{L^q}^q.
\end{equation}
\end{lem}

\noindent
\begin{prof}
The holomorphicity on $D(0,\epsilon)^{\ast}$ is evident. Let us 
show \eqref{eq4,3}. 

Thanks to \eqref{eq4,1}, we have
\begin{equation}\label{eq4,4}
U \big( \hd(b,0) - \mu \big)^{-1} \textup{P}_\perp = 
-U p \begin{pmatrix}
\mu^{-1} & 0\\
0 & 0
\end{pmatrix}.
\end{equation}  
As in \cite[Proof of Lemma 2.4]{rage}, it can proved that 
$U p \in \sqq \big( L^{2}(\br^2) \big)$ with
\begin{equation}\label{eq4,5}
\Vert U p \Vert_\sqq^q \leq
\frac{b_0}{2\pi} e^{2\textup{osc} \hspace{0.5mm} 
\tilde{\varphi}} \Vert U \Vert_{L^q}^q.
\end{equation}
Then, bound \eqref{eq4,3} follows by combining 
\eqref{eq4,4} and \eqref{eq4,5}.
\end{prof} 

\noindent
The next result concerns the second term of the r.h.s. 
of \eqref{eq4,2}.

\begin{lem}\label{l4,2} 
Under the assumptions of Lemma \ref{l4,1}, 
the operator-valued function 
$$
\bc \setminus [\zeta,+\infty) \ni \mu \longmapsto 
U \big( \hd(b,0) - \mu \big)^{-1} \textup{Q}_\perp
$$
is holomorphic with values in $\sqq \big( L^{2}(\br^2) \big)$. 
Moreover,
\begin{equation}\label{eq4,6}
\left\Vert U \big( \hd(b,0) - \mu \big)^{-1} \textup{Q}_\perp 
\right\Vert_\sqq^q \leq C \Vert U \Vert_{L^{q}}^{q} 
\left( 1 + \frac{\vert \mu + 1 \vert}{{\rm dist} \big( \mu,[\zeta,+\infty) \big)} 
\right)^q,
\end{equation}
where $C = C(q)$ is a constant depending on $q$.
\end{lem}

\noindent
\begin{prof}
For $\mu$ lying in the resolvent set of $\hd(b,0)$,
\eqref{eq1,1} implies that we have
\begin{equation}\label{eq4,7}
\left( \hd(b,0) - \mu \right)^{-1} \textup{Q}_{\perp} =
\left( \begin{smallmatrix}
   (\hd^- - \mu)^{-1} (I - p) & 0\\
   0 & (\hd^+ - \mu)^{-1}
\end{smallmatrix} \right) = (\hd^- - \mu)^{-1} (I - p) \oplus (\hd^+ - \mu)^{-1}.
\end{equation}
Then, $\bc \setminus [\zeta,+\infty) \ni \mu \longmapsto 
(\hd^- - \mu)^{-1} (I - p) \oplus (\hd^+ - \mu)^{-1}$ is well defined 
and analytic since $\bc \setminus [\zeta,+\infty)$ is included in the 
resolvent set of $\hd^-$ defined on $(I - p) Dom(\hd^-)$ 
and $\hd^+$ defined on $Dom(\hd^+)$. Thus,
$\mu \longmapsto U \big( \hd(b,0) - \mu \big)^{-1} \textup{Q}_\perp$
is holomorphic on $\bc \setminus [\zeta,+\infty)$.

Now, let us prove estimate \eqref{eq4,6}. In what follows below, 
constants are generic. Namely changing from a relation to another. 
First, let us show that \eqref{eq4,6} holds if $q$ is even. 

Identity \eqref{eq4,7} implies that we have
\begin{equation}\label{eq4,8}
\left\Vert U \left( \hd(b,0) - \mu \right)^{-1} \textup{Q}_{\perp} 
\right\Vert_\sqq^q \leq 
\left\Vert U \left( \hd^- - \mu \right)^{-1} (I - p) \right\Vert_\sqq^q
+ \left\Vert U \left( \hd^+ - \mu \right)^{-1} \right\Vert_\sqq^q.
\end{equation}
Let us focus on the first term of the r.h.s. of \eqref{eq4,8}. We have
\begin{equation}\label{eq4,9}
\left\Vert U (\hd^- - \mu)^{-1} (I - p) \right\Vert_\sqq^q 
\leq  \left\Vert U (\hd^- + 1)^{-1} \right\Vert_\sqq^q
\left\Vert (\hd^- + 1)(\hd^- - \mu)^{-1} (I - p) \right\Vert^q.
\end{equation}
By the Spectral mapping theorem, we have
\begin{equation}\label{eq4,10}
\left\Vert (\hd^- + 1)(\hd^- - \mu)^{-1} (I - p) \right\Vert^q \leq 
\textup{sup}_{s \in [\zeta,+\infty)}^q \left\vert \frac{s + 1}{s - \mu} 
\right\vert.
\end{equation}
Using the resolvent equation, the boundedness of the 
magnetic field $b$, and applying the diamagnetic inequality 
\big(see \cite[Theorem 2.3]{avr} and \cite[Theorem 2.13]{sim}\big)
which is only valid when $q$ is even, we get
\begin{equation}\label{eq4,11}
\begin{split}
\left\Vert U \bigl( \hd^- + 1 \bigr)^{-1} \right\Vert^q_\sqq 
& \leq \left\Vert I + (\hd^- + 1)^{-1}b \right\Vert^q
\left\Vert U \big( (-i\nabla - \textbf{A})^{2} + 1 \big)^{-1} 
\right\Vert^q_\sqq \\
& \leq C \left\Vert U (-\Delta + 1)^{-1} \right\Vert^q_\sqq.
\end{split}
\end{equation}
By the standard criterion \cite[Theorem 4.1]{sim}, we have
\begin{equation}\label{eq4,12}
\left\Vert U (-\Delta + 1 \vert)^{-1} \right\Vert^q_\sqq 
\leq C \Vert U \Vert_{L^q}^q 
\left\Vert \Bigl( \vert \cdot \vert^{2} + 1 \Bigr)^{-1} \right\Vert_{L^q}^q.
\end{equation}
Combining \eqref{eq4,9}, \eqref{eq4,10}, \eqref{eq4,11} and
\eqref{eq4,12}, we obtain
\begin{equation}\label{eq4,13}
\begin{split}
\left\Vert U \bigl( \hd^- - \mu \bigr)^{-1} (I - p) \right\Vert^q_\sqq
& \leq C \Vert U \Vert_{L^{q}}^{q} 
\textup{sup}_{s \in [\zeta,+\infty)}^q \left\vert \frac{s + 1}{s - \mu} 
\right\vert \\
& \leq C \Vert U \Vert_{L^{q}}^{q} \left( 1 + \frac{\vert \mu + 1 \vert}
{{\rm dist} \big( \mu,[\zeta,+\infty) \big)} \right)^q.
\end{split}
\end{equation}
Arguing similarly, it can be proved that we have
\begin{equation}\label{eq4,14}
\left\Vert U \bigl( \hd^+ - \mu \bigr)^{-1} \right\Vert^q_\sqq
\leq C \Vert U \Vert_{L^{q}}^{q} 
\left( 1 + \frac{\vert \mu + 1 \vert}
{{\rm dist} \big( \mu,[\zeta,+\infty) \big)} \right)^q.
\end{equation}
Hence, for $q$ even, \eqref{eq4,6} holds by putting together
\eqref{eq4,8}, \eqref{eq4,13} and \eqref{eq4,14}.

To prove that \eqref{eq4,6} holds for any $q \geq 2$, we use 
interpolation method. If $q$ verifies $q > 2$, then clearly there 
exists even integers $q_{0} < q_{1}$ satisfying $q \in (q_{0},q_{1})$ 
and $q_{0} \geq 2$. Consider $\gamma \in (0,1)$ such that
$\frac{1}{q} = \frac{1-\gamma}{q_{0}} + \frac{\gamma}{q_{1}}$ and 
introduce the operator
$$
L^{q_i} \big( \br^2 \big) \ni U \overset{T}{\longmapsto}  
U \big( \hd(b,0) - \mu \big)^{-1} \textup{Q}_\perp \in \mathcal{S}_{q_{i}},
\qquad i = 0, 1.
$$
Denote $C_{i} = C(q_{i})$ the constant appearing in \eqref{eq4,6} 
for $i = 0$, $1$, and define
$$
C(\mu,q_{i}) := C_i^{\frac{1}{q_i}} 
\left( 1 + \frac{\vert \mu + 1 \vert}{{\rm dist} \big( \mu,[\zeta,+\infty) \big)} \right).
$$
By \eqref{eq4,6},
$\Vert T \Vert \leq C(\mu,q_{i})$ for $i = 0$, $1$. Using the 
Riesz-Thorin Theorem \big(see for instance \cite[Sub. 5 of Chap. 6]{fol}, 
\cite{rie}, \cite{tho}, \cite[Chap. 2]{lun}\big),
we interpolate between $q_{0}$ and $q_{1}$ to obtain the extension 
$T : L^{q}(\br^2) \longrightarrow \sqq$, with
$$
\Vert T \Vert \leq C(\mu,q_{0})^{1-\gamma} C(\gamma,q_{1})^{\gamma} 
\leq C(q)^{\frac{1}{q}} \left( 1 + \frac{\vert \mu + 1 \vert}
{{\rm dist} \big( \mu,[\zeta,+\infty) \big)} \right).
$$
Then, in particular, we have for any $U \in L^q(\br^2)$
$$
\Vert T(U) \Vert_\sqq \leq 
C(q)^{\frac{1}{q}} \left( 1 + \frac{\vert \mu + 1 \vert}
{{\rm dist} \big( \mu,[\zeta,+\infty) \big)} \right)
\Vert U \Vert_{L^q},
$$
or equivalently bound \eqref{eq4,6}, and the lemma follows.
\end{prof}

Now, observe that Assumption (A1) given by \eqref{eq1,12}, 
implies that there exists a bounded operator $\mathscr{M}$ 
such that $\vert V \vert^\frac{1}{2} = \mathscr{M} F^\frac{1}{2}$
with $F$ satisfying \eqref{eq1,12}. Consequently, since the 
radius $\epsilon$ of $D(0,\epsilon)^\ast$ satisfies 
$0 < \epsilon < \zeta$ and $V$ is bounded, then by combining 
identities \eqref{eq4,1} and \eqref{eq4,2} with Lemmas \ref{l4,1} 
and \ref{l4,2}, we obtain the following:

\begin{lem}\label{l4,3}
The operator-valued function
$$
D(0,\epsilon)^{\ast} \ni \mu \longmapsto 
\mathcal{T}_{V,\perp}(\mu) 
:= \Tilde{J} \vert V \vert^{\frac{1}{2}} \big( \hd(b,0) - 
\mu \big)^{-1} \vert V \vert^{\frac{1}{2}}
$$ 
is holomorphic with values in $\sqq \big( L^{2}(\br^2) \big)$,
where $\Tilde{J}$ is defined by the polar decomposition 
$V = \Tilde{J} \vert V \vert$ of $V$.
\end{lem}

\subsection{Reduction of the problem}

We reduce the study of the discrete eigenvalues of $\hd(b,V)$ near 
zero, to that of the zeros of a holomorphic function in a 
punctured neighbourhood of zero.

In what follows below, the regularized determinant 
$\textup{det}_{\lceil q \rceil} (\bullet)$ is defined in Appendix 
A $1$. Similarly to Lemma \ref{l4,3}, it can be shown that 
$V \big( \hd(b,0) - \cdot \big)^{-1}$ is holomorphic 
in $D(0,\epsilon)^{\ast}$ with values in 
$\sqq \big( L^{2}(\br^2) \big)$. Then,
$\textup{det}_{\lceil q \rceil} \Bigl( I + V \big( \hd(b,0) - 
\mu \big) \Bigr)$ is well defined by \eqref{eq3,2}. It is well known, 
see for instance \cite[Chap. 9]{sim}, that  
\begin{equation}\label{eq4,15}
\mu \in \sigma_{\textup{\textbf{disc}}} \big( \hd (b,V) \big) 
\Leftrightarrow f_{q}(\mu) := 
\textup{det}_{\lceil q \rceil} \Bigl( I + V \big( \hd(b,0) - 
\mu \big)^{-1} \Bigr) = 0.
\end{equation}
Since $V \big( \hd(b,0) - \cdot \big)^{-1}$ is holomorphic in  
$D(0,\epsilon)^{\ast}$, then so is the function $f_q$ 
by Property \textbf{d)} of Appendix A $1$. Moreover, the 
algebraic multiplicity of $\mu$ as a discrete eigenvalue of 
$\hd(b,V)$ is equal to its order as a zero of $f_q$.

In the next proposition, the quantity $Ind_{\mathscr{C}}(\bullet)$
appearing in the r.h.s. of \eqref{eq4,16} is defined in Appendix 
A $2$.

\begin{prop}\label{p4,1} 
Let $\mathcal{T}_{V,\perp} (\mu)$ be the operator defined in 
Lemma \ref{l4,3}. Then, the following assertions are equivalent:
\begin{itemize}
\item[(i)] $\mu \in D(0,\epsilon)^{\ast}$ is a 
discrete eigenvalue of $\hd(b,V)$,
\item[(ii)] $\textup{det}_{\lceil q \rceil} 
\big( I + \mathcal{T}_{V,\perp} (\mu) \big)= 0$,
\item[(iii)] $-1$ is an eigenvalue of $\mathcal{T}_{V,\perp}(\mu)$.

\medskip

\noindent
Moreover,
\begin{equation}\label{eq4,16}
\textup{mult}(\mu) = 
Ind_{\mathscr{C}} \hspace{0.5mm} \big( I + \mathcal{T}_{V,\perp}
(\cdot) \big),
\end{equation}
where $\mathscr{C}$ is a small contour positively oriented,
containing $\mu$ as the unique discrete eigenvalue of $\hd(b,V)$.
\end{itemize}
\end{prop}

\noindent
\begin{prof}
The equivalence between assertions (i) and (ii) follows 
obviously from \eqref{eq4,15} and the equality
$$
\textup{det}_{\lceil q \rceil} \Bigl( I + V \big( \hd(b,0) - 
\mu \big) \Bigr)
= \textup{det}_{\lceil q \rceil} \Bigl( I + \Tilde{J} 
\vert V \vert^\frac{1}{2}
\big( \hd(b,0) - \mu \big)^{-1} \vert V \vert^\frac{1}{2} \Bigr),
$$
thanks to Property \textbf{b)} of Appendix A $1$.

The equivalence between assertions (ii) and (iii) 
follows from Property \textbf{c)} of Appendix A $1$. 

We prove now \eqref{eq4,16}. To this end, we consider the function 
$f_q$ defined by \eqref{eq4,15}. According to the comment just 
after \eqref{eq4,15}, we have
\begin{equation}\label{eq4,17}
\textup{mult}(\mu) = ind_{\mathscr{C}} f_q,
\end{equation}
where the r.h.s. of \eqref{eq4,17} is the index defined by \eqref{eqa,1} 
of the holomorphic function $f_q$ with respect to the contour 
$\mathscr{C}$. We thus get \eqref{eq4,16} easily from the equality
$$
ind_{\mathscr{C}} f_q = Ind_{\mathscr{C}} \hspace{0.5mm} 
\Big( I + \mathcal{T}_{V,\perp}(\cdot) \Big),
$$
\big(see for instance \cite[(2.6)]{bo} for more details\big).
\end{prof}

\subsection{Decomposition of the weighted resolvent}

Our goal in this section is to split the weighted resolvent 
$\mathcal{T}_{V,\perp}(\mu) := \Tilde{J} 
\vert V \vert^{\frac{1}{2}} \big( \hd(b,0) - \mu \big)^{-1} 
\vert V \vert^{\frac{1}{2}}$ into two parts which are 
respectively meromorphic and holomorphic in a neighbourhood of 
zero. The potential $V$ is assumed to satisfy Assumption (A1).

The next proposition is a direct consequence of identities 
\eqref{eq4,1}, \eqref{eq4,2} and Lemma \ref{l4,3}.

\begin{prop}\label{p4,2} 
For $\mu \in D(0,\epsilon)^{\ast}$, we have
\begin{equation}\label{eq4,18}
\mathcal{T}_{V,\perp}(\mu) = 
-\frac{\Tilde{J}}{\mu} \vert V \vert^{\frac{1}{2}}
\begin{pmatrix}
p & 0 \\
0 & 0
\end{pmatrix}
\vert V \vert^{\frac{1}{2}} + \mathscr{A}_{\perp}(\mu),
\end{equation}
where the operator
$\mathscr{A}_{\perp}(\mu) := \Tilde{J} \vert V \vert^{\frac{1}{2}} 
\big( \hd(b,0) - \mu \big)^{-1} \textup{Q}_\perp
\vert V \vert^{\frac{1}{2}} \in \sqq \big( L^{2}(\br^2) \big)$
is holomorphic in the open disk 
$D(0,\epsilon) := D(0,\epsilon)^{\ast} \cup \lbrace 0 \rbrace$.
\end{prop} 

\begin{rem}\label{r4,1}
\end{rem}

\begin{itemize}
\item[(i)] Thanks to Lemma \ref{l4,1},
$\vert V \vert^{\frac{1}{2}}
\begin{pmatrix}
p & 0 \\
0 & 0
\end{pmatrix}
\vert V \vert^{\frac{1}{2}}$ is a compact operator. Then,
for any $r > 0$, we have
\begin{equation}\label{eq4,19}
\textup{Tr} \hspace{0.4mm} \one_{(r,\infty)} 
\left( \vert V \vert^{\frac{1}{2}}
\begin{pmatrix}
p & 0 \\
0 & 0
\end{pmatrix}
\vert V \vert^{\frac{1}{2}} \right) 
= \textup{Tr} \hspace{0.4mm} \one_{(r,\infty)} 
\left( \begin{pmatrix}
p & 0 \\
0 & 0
\end{pmatrix}
\vert V \vert \begin{pmatrix}
p & 0 \\
0 & 0
\end{pmatrix}\right)
= \textup{Tr} \hspace{0.4mm} \one_{(r,\infty)} 
\left( p \vert V \vert_{11} p  \right),
\end{equation}
where we recall that the $\vert V \vert_{\ell k}$, $1 \leq 
\ell,k \leq 2$, are the coefficients of the matrix 
$\vert V \vert$.

\item[(ii)] If $V$ satisfies Assumption (A2) given by 
\eqref{eq2,7}, then Proposition \ref{p4,2} holds with $\Tilde{J}$ 
replaced by $J\eta$, $J := sign(W)$, and $\vert V \vert_{11}$
replaced by $\vert W \vert_{11}$ in \eqref{eq4,19}.
\end{itemize}

\section{Proof of Theorem \ref{t2,1}: Upper bound}\label{s5}

The proof will be divided into two steps.

\subsection{A preliminary result}

Let
\begin{equation}\label{eq:n}
N \big( \hd(b,V) \big) := \big\lbrace \langle \hd(b,V)f,f \rangle : 
f \in Dom \big( \hd(b,V) \big), \Vert f \Vert_{L^{2}} = 1 \big\rbrace,
\end{equation}
denote the numerical range of the operator $\hd(b,V)$. The inclusion 
$\sigma \big( \hd(b,V) \big) \subseteq 
\overline{N \big( \hd(b,V) \big)}$ is well known
(see for instance \cite[Lemma 9.3.14]{dav}\big).

\begin{prop}\label{p5,1} 
There exists $r_{0} > 0$ such that for any 
$0 < r < \vert \mu \vert < r_{0}$, we have:

\begin{itemize} 
\item[(i)] $\mu$ is a discrete eigenvalue of 
$\hd(b,V)$ near zero if and only if $\mu$ is a zero of 
\begin{equation}\label{eq5,1}
\mathscr{D}_{\perp}(\mu,r) := \det \big( I + 
\mathscr{K}_{\perp}(\mu,r) \big),
\end{equation}
with $\mathscr{K}_{\perp}(\mu,r)$ a finite-rank operator 
analytic with respect to $\mu$ and satisfying
$$
\textup{rank} \hspace{0.6mm} \mathscr{K}_{\perp}(\mu,r) = 
\mathcal{O} \Big( \textup{Tr} \hspace{0.4mm} 
\one_{(r,\infty)} \big( p \vert V \vert_{11} p 
\big) + 1 \Big), \hspace{0.3cm} 
\left\Vert \mathscr{K}_{\perp}(\mu,r) \right\Vert = 
\mathcal{O} \big( r^{-1} \big),
$$
where the $\mathcal{O}$'s are uniform with respect to 
$r$,  $\mu$.

\item[(ii)] Furthermore, if $\mu$ is a discrete 
eigenvalue of $\hd(b,V)$ near zero, then,
\begin{equation}\label{eq5,2}
\textup{mult}(\mu) = Ind_\mathscr{C} \hspace{0.5mm} 
\left( I + \mathscr{K}_{\perp}(\cdot,r) \right) = 
\textup{m}(\mu),
\end{equation}
$\mathscr{C}$ being chosen as in \eqref{eq4,16}, and 
$\textup{m}(\mu)$ being the multiplicity of $\mu$ as a zero of 
$\mathscr{D}_{\perp}(\cdot,r)$.

\item[(iii)] If $\mu$ verifies $\textup{dist} 
\big( \mu,\overline{N \big( \hd(b,V) \big)} \big) > \varsigma > 0$,
$\varsigma = \mathcal{O}(1)$, then $I + \mathscr{K}_{\perp}(\mu,r)$ is 
invertible and satisfies
$
\left\Vert \big( I + \mathscr{K}_{\perp}(\mu,r) \big)^{-1} 
\right\Vert = \mathcal{O} \left( \varsigma^{-1} \right)
$, where the $\mathcal{O}$ is uniform with respect to 
$r$,  $\mu$ and $\varsigma$.
\end{itemize}
\end{prop}

\noindent
\begin{prof}
\textit{(i)-(ii):} By Proposition \ref{p4,2}, $\mu \mapsto 
\mathscr{A}_\perp (\mu)$ is holomorphic near zero with values in 
$\sqq \big( L^2(\br^2) \big)$. Then, for $r_{0}$ small enough, 
there exists a finite-rank operator $\mathscr{A}_{0,\perp}$ 
independent of $\mu$, and an operator $\tilde{\mathscr{A}}_\perp 
(\mu) \in \sqq \big( L^2(\br^2) \big)$ holomorphic near zero satisfying 
$\Vert \tilde{\mathscr{A}}_{\perp} (\mu) \Vert < \frac{1}{4}$
for $0 < \vert \mu \vert \leq r_{0}$, and
$$
\mathscr{A}_\perp(\mu) = \mathscr{A}_{0,\perp} + 
\tilde{\mathscr{A}}_\perp(\mu).
$$ 
Set $\mathscr{B}_\perp := \vert V \vert^{\frac{1}{2}}
\begin{pmatrix}
p & 0 \\
0 & 0
\end{pmatrix}
\vert V \vert^{\frac{1}{2}}$ and write
\begin{equation}\label{eq5,3}
\mathscr{B}_\perp = \mathscr{B}_\perp \one_{[0,\frac{1}{2}r]} 
(\mathscr{B}_\perp) + \mathscr{B}_\perp 
\one_{(\frac{1}{2}r,\infty)} (\mathscr{B}_\perp).
\end{equation}
It is easy to check that 
$\left\Vert -\frac{\Tilde{J}}{\mu} \mathscr{B}_\perp 
\one_{[0,\frac{1}{2}r]} (\mathscr{B}_\perp) 
+ \tilde{\mathscr{A}}_\perp(\mu) 
\right\Vert < \frac{3}{4}$ for $0 < r < \vert \mu \vert < r_{0}$. 
Consequently, we have
\begin{equation}\label{eq5,4}
\big( I + \mathcal{T}_{V,\perp}(\mu) \big) = 
\big( I + \mathscr{K}_\perp(\mu,r) \big) 
\left( I - \frac{\Tilde{J}}{\mu} \mathscr{B}_\perp 
\one_{[0,\frac{1}{2}r]} (\mathscr{B}_\perp) 
+ \tilde{\mathscr{A}}_\perp(\mu) \right),
\end{equation}
where $\mathscr{K}_\perp(\mu,r)$ is defined by
$$
\small{\mathscr{K}_\perp(\mu,r) := 
\left( -\frac{\Tilde{J}}{\mu} \mathscr{B}_\perp 
\one_{(\frac{1}{2}r,\infty)} (\mathscr{B}_\perp) 
+ \mathscr{A}_{0,\perp} \right) 
\left( I - \frac{\Tilde{J}}{\mu} \mathscr{B}_\perp 
\one_{[0,\frac{1}{2}r]} (\mathscr{B}_\perp) 
+ \tilde{\mathscr{A}}_\perp(\mu) \right)^{-1}}.
$$
It is a finite-rank operator of order
$$
\small{\mathcal{O} \Big( \textup{Tr} \hspace{0.4mm} 
\one_{(\frac{1}{2}r,\infty)} (\mathscr{B}_\perp) 
+ 1 \Big) = \mathcal{O} \Big( \textup{Tr} \hspace{0.4mm} 
\one_{(r,\infty)} \big( p \vert V \vert_{11} p 
\big) + 1 \Big)},
$$
taking into account \eqref{eq4,19}. Moreover, its norm 
is of order $\mathcal{O} \big( \vert \mu \vert^{-1} \big) = 
\mathcal{O} \big( r^{-1} \big)$. Since we have
$\Vert - \frac{\Tilde{J}}{\mu} \mathscr{B}_\perp 
\one_{[0,\frac{1}{2}r]} (\mathscr{B}_\perp) 
+ \tilde{\mathscr{A}}_\perp (\mu)\Vert < 1$ 
for $0 < r < \vert \mu \vert < r_{0}$, then 
\cite[Theorem 4.4.3]{goh} implies that
$$
\small{Ind_\mathscr{C} \hspace{0.5mm} \left(I - 
\frac{\Tilde{J}}{\mu} \mathscr{B}_\perp \one_{[0,\frac{1}{2}r]} 
(\mathscr{B}_\perp) + \tilde{\mathscr{A}}_\perp (\mu) \right) = 0}.
$$ 
Therefore, Property \eqref{eqa,3} applied to \eqref{eq5,4} together 
with Proposition \ref{p4,1} give \eqref{eq5,2}. Furthermore, it follows 
that $\mu$ is a discrete eigenvalue of $\hd(b,V)$ near zero if and 
only if $\mu$ is a zero of $\mathscr{D}_\perp(\cdot,r)$.

\textit{(iii):} Equality \eqref{eq5,4} implies that we have
\begin{equation}\label{eq5,5}
I + \mathscr{K}_\perp(\mu,r) = \big( I + \mathcal{T}_{V,\perp}(\mu) 
\big) \left( I - \frac{\Tilde{J}}{\mu} \mathscr{B}_\perp 
\one_{[0,\frac{1}{2}r]} (\mathscr{B}_\perp) 
+ \tilde{\mathscr{A}}_\perp(\mu) \right)^{-1},
\end{equation}
for $0 < r < \vert \mu \vert < r_{0}$. From the resolvent equation,
it is easy to deduce that
$$
\small{\left( I + \Tilde{J} \vert V \vert^{1/2} 
\big( \hd(b,0) - \mu \big)^{-1} \vert V \vert^{1/2} \right) 
\left( I - \Tilde{J} \vert V \vert^{1/2} \big( \hd(b,V) - \mu \big)^{-1} 
\vert V \vert^{1/2} \right) = I}.
$$
Then, if $\mu$ belongs to resolvent set of $\hd(b,V)$, we have
$$
\big( I + \mathcal{T}_{V,\perp}(\mu) \big)^{-1} 
= I - \Tilde{J} \vert V \vert^{1/2} \big( \hd(b,V) - \mu \big)^{-1} 
\vert V \vert^{1/2}.
$$
This together with \eqref{eq5,5} imply that $I + \mathscr{K}_\perp(\mu,r)$ 
is invertible for $0 < r < \vert \mu \vert < r_{0}$. So, using
\cite[Lemma 9.3.14]{dav}, we obtain
\begin{align*}
\left\Vert \big( I + \mathscr{K}_\perp(\mu,r) \big)^{-1} \right\Vert 
& = \mathcal{O} \Big( 1 + \left\Vert \vert V \vert^{1/2} 
\big( \hd(b,V) - \mu \big)^{-1} \vert V \vert^{1/2} \right\Vert \Big) \\
&= \mathcal{O} \Big( 1 + \textup{dist} 
\big( \mu,\overline{N \big( \hd(b,V) \big)} \big)^{-1} \Big)
= \mathcal{O} \left( \varsigma^{-1} \right),
\end{align*}
whenever
$\textup{dist} \big( \mu,\overline{N \big( \hd(b,V) \big)} \big) 
> \varsigma > 0$, $\varsigma = \mathcal{O}(1)$,
and the proof is complete.
\end{prof}

\subsection{Back to the proof of Theorem \ref{t2,1}}

Proposition \ref{p5,1} implies that for $0 < r < \vert \mu \vert < r_{0}$, 
we have
\begin{equation}\label{eq5,6}
\begin{aligned}
\mathscr{D}_\perp(\mu,r) & = \prod_{j=1}^{\mathcal{O} 
\big( \textup{Tr} \hspace{0.4mm} 
\one_{(r,\infty)} (p \vert V \vert_{11} p) 
+ 1 \big)} \big( 1 + \lambda_{j}(\mu,r) \big) \\
& = \mathcal{O}(1) \hspace{0.5mm} \textup{exp} \hspace{0.5mm} 
\Big( \mathcal{O} \big( \textup{Tr} \hspace{0.4mm} 
\one_{(r,\infty)} 
\big( p \vert V \vert_{11} p \big) + 1 \big) 
\vert \ln r \vert \Big),
\end{aligned}
\end{equation}
$\lambda_{j}(\mu,r)$ being the eigenvalues of 
$\mathscr{K}_\perp := \mathscr{K}_\perp(\mu,r)$ 
satisfying
$\vert \lambda_{j}(\mu,r) \vert = 
\mathcal{O} \left( r^{-1} \right)$. Let 
$\mu \in D(0,\epsilon)^\ast$ satisfy $\textup{dist} 
\big( \mu,\overline{N \big( \hd(b,V) \big)} \big) > \varsigma > 0$ 
and $0 < r < \vert \mu \vert < r_{0}$. Then, we have
$$
\mathscr{D}_\perp(\mu,r)^{-1} = 
\det \big( I + \mathscr{K}_\perp \big)^{-1} = 
\det \big( I - \mathscr{K}_\perp(I+\mathscr{K}_\perp)^{-1} \big).
$$
Similarly to \eqref{eq5,6}, it can be shown that
\begin{equation}\label{eq5,7}
\vert \mathscr{D}_\perp(\mu,r) \vert \geq 
C \hspace{0.5mm} \textup{exp} \hspace{0.5mm} 
\Big( - C \big( \textup{Tr} \hspace{0.4mm} 
\one_{(r,\infty)} 
\big( p \vert V \vert_{11} p \big) + 1 \big) 
\big( \vert \textup{ln} \hspace{0.5mm} \varsigma \vert 
+ \vert \textup{ln} \hspace{0.5mm} r \vert \big) \Big),
\end{equation}
so that for $\frac{1}{4}r < \varsigma < 2r$, $0 < r \ll 1$, we obtain
\begin{equation}\label{eq5,70}
- \ln \, \vert \mathscr{D}_\perp(\mu,r) \vert \leq
C \, \textup{Tr} \hspace{0.4mm} \one_{(r,\infty)} 
\big( p \vert V \vert_{11} p \big) \vert \textup{ln} \hspace{0.5mm} r \vert
+ \mathcal{O}(1).
\end{equation}
Consider the discrete eigenvalues
$\mu \in \big\lbrace r < \vert \mu \vert < 2r \big\rbrace 
\subset D(0,\zeta)^{\ast}$ with $r > 0$ such that $r < \Vert V \Vert < 
\frac{3}{2}r$. Thanks to their discontinuous distribution, there exists 
a simply connected sub-domain $\Delta$ of $\big\lbrace r < \vert \mu \vert < 2r 
\big\rbrace$ containing all the eigenvalues $\mu$ and such that 
$\sigma_{\text{\textbf{disc}}} \big( \hd(b,V) \big) \cap
\partial\Delta = \varnothing$. Note that in view of the definition \eqref{eq:n}
of the numerical range $N \big( \hd(b,V) \big)$ of the operator $\hd(b,V)$, we
have 
\begin{equation}
N \big( \hd(b,V) \big) \subseteq \big\lbrace \mu \in \bc : \vert \Im(\mu) 
\vert \leq \Vert V \Vert \big\rbrace.
\end{equation}
Then, Theorem \ref{t2,1} holds by applying the Jensen Lemma \ref{la,1} with 
the function  $g(\mu) := \mathscr{D}_\perp(r \mu,r)$, $\mu \in \Delta / r$, 
with some $\mu_{0} \in \Delta/r$ satisfying $\textup{dist} \big( r\mu_0,
\overline{N \big( \hd(b,V) \big)} \big) \geq \varsigma > \frac{1}{4}r$, 
$\varsigma < 2r$, and by using \eqref{eq5,6} and \eqref{eq5,70}.

\section{Proof of Theorem \ref{t2,2}: Localisation and asymptotics expansions}\label{s6}

First, we have to rephrase Proposition \ref{p4,1} with respect to the 
characteristic value terminology (see Definition \ref{d1}).

\begin{prop}\label{p6,1} 
For $\mu \in D(0,\epsilon)^{\ast}$, the following assertions 
are equivalent:

\begin{itemize}
\item[(i)] $\mu$ is a discrete eigenvalue
of $\hd(b,V)$,

\item[(ii)] $\mu$ is a characteristic 
value of $I + \mathcal{T}_{V,\perp}(\cdot)$.

\smallskip

\noindent
Moreover, the multiplicity of $\mu$ as a discrete eigenvalue
coincides with its multiplicity as a characteristic value 
defined by \eqref{eqa,6}.
\end{itemize}
\end{prop}

\noindent
\textit{Proof of assertion (i) of Theorem \ref{t2,2}:} From 
Proposition \ref{p6,1} and according to (ii) of Remark \ref{r4,1}, 
we reduce the investigation of the discrete eigenvalues 
$\mu \in D(0,\epsilon)^{\ast}$ to that of the characteristic 
values of 
$$
I + \mathcal{T}_{V,\perp}(\mu) = I - \eta 
\frac{\mathcal{A}_\perp(\mu)}{\mu},
$$ 
the operator $\mathcal{A}_\perp(\mu)$ being defined by 
\eqref{eq2,9}. In particular,
$\pm \mathcal{A}_\perp(0) = \vert W \vert^{\frac{1}{2}}
\begin{pmatrix}
p & 0 \\
0 & 0
\end{pmatrix}
\vert W \vert^{\frac{1}{2}}$
for $J = \pm$. 
Thus, assertion (i) of Theorem \ref{t2,1} holds by (i) and (ii) of 
Lemma \ref{la,2} with $z = \pm \mu/\eta$. To be more precise, near zero, 
the discrete eigenvalues $\mu$ verify for any $\delta > 0$
\begin{equation}\label{eq6,1}
\pm \Re \left( \frac{\mu}{\eta} \right) \geq 0, 
\hspace{0.7cm}  \mu \in \pm \eta \hspace{0.05cm} 
\overline{\Gamma^{\delta}(r,r_{0})},
\end{equation}
the sector $\Gamma^{\delta}(r,r_{0})$ being defined by 
\eqref{eq2,5}.

\medskip

\noindent
\textit{Proof of assertion (ii) of Theorem \ref{t2,2}:} The above 
proof of (i) of Theorem \ref{t2,2} together with Proposition 
\ref{p6,1}, imply that the discrete eigenvalues $\mu$ near zero are 
the characteristic values 
$\mu \in \mathcal{Z} \big( D(0,\epsilon)^{\ast} \big)$ of
$I + \mathcal{T}_{V,\perp}(\cdot)$ concentrated
in the sectors $\big\lbrace \mu \in D(0,\epsilon)^{\ast} : 
\pm \mu/\eta \in \Gamma^{\delta}(r,r_{0}) \big\rbrace$,
for any $\delta > 0$. In particular, we obtain
\begin{equation}\label{eq6,2}
\begin{split}
\displaystyle \sum_{\substack{\mu \hspace{0.5mm} \in 
\hspace{0.5mm} \sigma_{\text{\textbf{disc}}} \big( \hd(b,V) 
\big) \\
r < \vert \mu \vert < r_{0}}} \textup{mult}(\mu) 
 = \displaystyle \sum_{\substack{\mu \in 
\mathcal{Z} \big( D(0,\epsilon)^{\ast} \big) \\ 
\pm \mu/\eta \hspace{0.5mm} \in \hspace{0.5mm} 
\Gamma^{\delta}(r,r_{0})}} \textup{mult}(\mu) + 
\mathcal{O}(1)  = \mathcal{N} 
\big( \Gamma^{\delta}(r,r_{0}) \big) + \mathcal{O}(1),
\end{split}
\end{equation}
the quantity $\mathcal{N}(\bullet)$ being defined by \eqref{eqa,7}. 
If $n(\bullet)$ is the quantity defined by \eqref{eqa,8} 
with $T(0) = \pm \mathcal{A}_\perp(0)$, then by using \eqref{eq4,19},
we get
\begin{equation}\label{eq6,3}
\begin{split}
n \big( [r,r_{0}] \big) = \textup{Tr} \hspace{0.4mm} 
\one_{(r,\infty)} \left( \vert W \vert^{\frac{1}{2}}
\begin{pmatrix}
p & 0 \\
0 & 0
\end{pmatrix}
\vert W \vert^{\frac{1}{2}} \right) 
+ \mathcal{O}(1) 
 = \textup{Tr} 
\hspace{0.6mm} \one_{(r,\infty)} 
\big( p \vert W \vert_{11} p \big) + \mathcal{O}(1).
\end{split}
\end{equation}
Thus, (ii) of Theorem \ref{t2,2} follows from (iii) of Lemma 
\ref{la,2} together with \eqref{eq6,2} and \eqref{eq6,3}.

\medskip

\noindent
\textit{Proof of assertion (iii) of Theorem \ref{t2,2}:} If we have 
$\Phi(r) = r^{-\gamma}$, or $\Phi(r) = 
\vert \ln r \vert^{\gamma}$, or $\Phi(r) = 
\big( \ln \vert \ln r \vert \big)^{-1} \vert \ln r \vert$
for some $\gamma > 0$, then it can be checked that
$$
\phi \big( r(1 \pm \nu) \big) = \phi (r) \big( 1 + o(1) 
+ \mathcal{O}(\nu) \big),
$$ 
for any $\nu > 0$ small enough. Then, (iii) of Theorem \ref{t2,2} holds by
(iv) of Lemma \ref{la,2} combined with \eqref{eq6,2} and \eqref{eq6,3}. 
This completes the proof.

\section{Discrete eigenvalues for the 3D problem}\label{s7}

\subsection{Preliminary results}

Define $P := p \otimes 1$, $Q := I - P$, and introduce the 
orthogonal projections in $L^2(\br^3)$
\begin{equation}\label{eq7,1}
\textup{P} := \begin{pmatrix} 
   P & 0 \\ 
   0 & 0 
\end{pmatrix}, \hspace{1cm} 
\textup{Q} := \textup{I} - \textup{P} = \begin{pmatrix} 
   Q & 0 \\
   0 & I 
\end{pmatrix}.
\end{equation}
For $z \notin \sigma \big( \htr(b,0) \big)$, on account of 
\eqref{eq1,1} with $n = 3$ and Proposition \ref{pr}, we have
\begin{equation}\label{eq7,2}
\big( \htr(b,0) - z \big)^{-1} \textup{P} = 
\begin{pmatrix}
p \otimes \mathscr{R}(z) & 0 \\
   0 & 0
\end{pmatrix},
\end{equation}
with the resolvent 
$\mathscr{R}(z) := \left( -\frac{d^2}{dX_\parallel^2} - z 
\right)^{-1}$ acting in $L^{2}(\br)$. Then, for any $z \in \bc 
\setminus [0,+\infty)$, we have
\begin{equation}\label{eq7,3}
\big( \htr(b,0) - z \big)^{-1} = 
\big( p \otimes \mathscr{R}(z) \big) 
\begin{pmatrix}
1 & 0\\
0 & 0
\end{pmatrix} + \big( \htr(b,0) - z \big)^{-1} \textup{Q}.
\end{equation}

\noindent
We have the following lemma:

\begin{lem}\label{l7,1}
For given $U \in L^q(\br^2)$ and $G \in L^q(\br)$, 
$q \in [2,+\infty)$, the operator-valued function 
$$
\bc \setminus [0,+\infty) \ni z \longmapsto U G 
\big( \htr(b,0) - z \big)^{-1} \textup{P} 
$$ 
is holomorphic with values in $\sqq \big( L^{2}(\br^3) \big)$.
Moreover,
\begin{equation}\label{eq7,4}
\left\Vert U G
\big( \htr(b,0) - z \big)^{-1} \textup{P} \right\Vert_\sqq^q \leq
C \frac{b_0 e^{2\textup{osc} \hspace{0.5mm} \tilde{\varphi}}}{2\pi}
\Vert U \Vert_{L^q}^q \Vert G \Vert_{L^{q}}^{q} 
\left( 1 + \frac{\vert z + 1 \vert}{{\rm dist} \big( z,[0,+\infty) \big)} 
\right)^q,
\end{equation}
where $C = C(q)$ is a constant depending only on $q$.
\end{lem}

\noindent
\begin{prof}
The holomorphicity on $\bc \setminus [0,+\infty)$ is trivial. 
Let us prove \eqref{eq7,4}. 

Thanks to \eqref{eq7,2}, we have
\begin{equation}\label{eq7,5}
U G \big( \htr(b,0) - z \big)^{-1} \textup{P} = 
\big( Up \otimes G\mathscr{R}(z) \big) 
\begin{pmatrix}
1 & 0\\
0 & 0
\end{pmatrix}.
\end{equation}  
As in \eqref{eq4,5}, we have
\begin{equation}\label{eq7,6}
\Vert U p \Vert_\sqq^q \leq
\frac{b_0 e^{2\textup{osc} \hspace{0.5mm} \tilde{\varphi}}}{2\pi}
\Vert U \Vert_{L^q}^q.
\end{equation}
From the estimate
\begin{equation}\label{eq7,7}
\left\Vert G \mathscr{R}(z) \right\Vert_\sqq^q 
\leq \left\Vert G \left( -\frac{d^2}{dX_\parallel^2} + 1 \right)^{-1} 
\right\Vert_\sqq^q
\left\Vert \left( -\frac{d^2}{dX_\parallel^2} + 1 \right)\mathscr{R}(z) 
\right\Vert^q,
\end{equation}
we obtain by the Spectral mapping theorem 
\begin{equation}\label{eq7,8}
\left\Vert \left( -\frac{d^2}{dX_\parallel^2} + 1 \right)\mathscr{R}(z) 
\right\Vert^q \leq \textup{sup}_{s \in [0,+\infty)}^q 
\left\vert \frac{s + 1}{s - z} \right\vert \leq
\left( 1 + \frac{\vert z + 1 \vert}{{\rm dist} \big( z,[0,+\infty) \big)} 
\right)^q,
\end{equation}
and by the standard criterion \cite[Theorem 4.1]{sim}, we obtain
\begin{equation}\label{eq7,9}
\left\Vert G \left( -\frac{d^2}{dX_\parallel^2} + 1 \right) \right\Vert^q_\sqq 
\leq C \Vert G \Vert_{L^q}^q 
\left\Vert \Bigl( \vert \cdot \vert^{2} + 1 \Bigr)^{-1} 
\right\Vert_{L^q}^q.
\end{equation}
Then, \eqref{eq7,4} follows by putting together 
\eqref{eq7,5}, \eqref{eq7,6}, \eqref{eq7,7}, \eqref{eq7,8} and
\eqref{eq7,9}.
\end{prof}

\noindent
The next lemma is just the analogue of Lemma \ref{l4,2} in 
dimension three. It can be proved in a similar way by taking 
into account the appropriate modifications. For this reason, 
to simplify our exposition, its proof will be omitted.

\begin{lem}\label{l7,2} 
For a given $g \in L^q(\br^3)$, $q \in [2,+\infty)$, the 
operator-valued function 
$$
\bc \setminus [\zeta,+\infty) \ni z \longmapsto 
g \big( \htr(b,0) - z \big)^{-1} \textup{Q}
$$
is holomorphic with values in $\sqq \big( L^{2}(\br^3) \big)$. 
Moreover,
\begin{equation}\label{eq7,10}
\left\Vert g \big( \htr(b,0) - z \big)^{-1} \textup{Q} 
\right\Vert_\sqq^q \leq C \Vert g \Vert_{L^{q}}^{q} 
\left( 1 + \frac{\vert z + 1 \vert}{{\rm dist} \big( z,[\zeta,+\infty) \big)} 
\right)^q,
\end{equation}
where $C = C(q)$ is a constant depending only on $q$.
\end{lem}

Throughout this article, we will use the following choice 
of the complex square root
\begin{equation}\label{eq2,16}
\bc \setminus (-\infty,0] \overset{\sqrt{\cdot}}{\longrightarrow} 
\bc_+.
\end{equation}
For $0 < \kappa < \sqrt{\zeta}$, let $D_\pm(0,\kappa^2)$ 
be the half-rings defined by \eqref{eq2,3b}. Put the change 
of variables $z = k^2$ and define the domains
\begin{equation}\label{eq2,18}
\mathcal{D}_\pm^\ast (\kappa) 
:= \big\lbrace k \in \bc_\pm : 0 < \vert k \vert < \kappa : 
\Re(k) > 0 \big\rbrace.
\end{equation}
Under the above considerations, $D_\pm(0,\kappa^2)$ can 
be parametrized by $z = z(k) := k^2$, with $k \in \mathcal{D}_\pm^\ast 
(\kappa)$ respectively (see Figure 6.1 below):

\begin{figure}[h]\label{fig 1}
\begin{center}

\vspace*{-1.5cm}

\hspace*{-1cm} \includegraphics[scale=0.6]{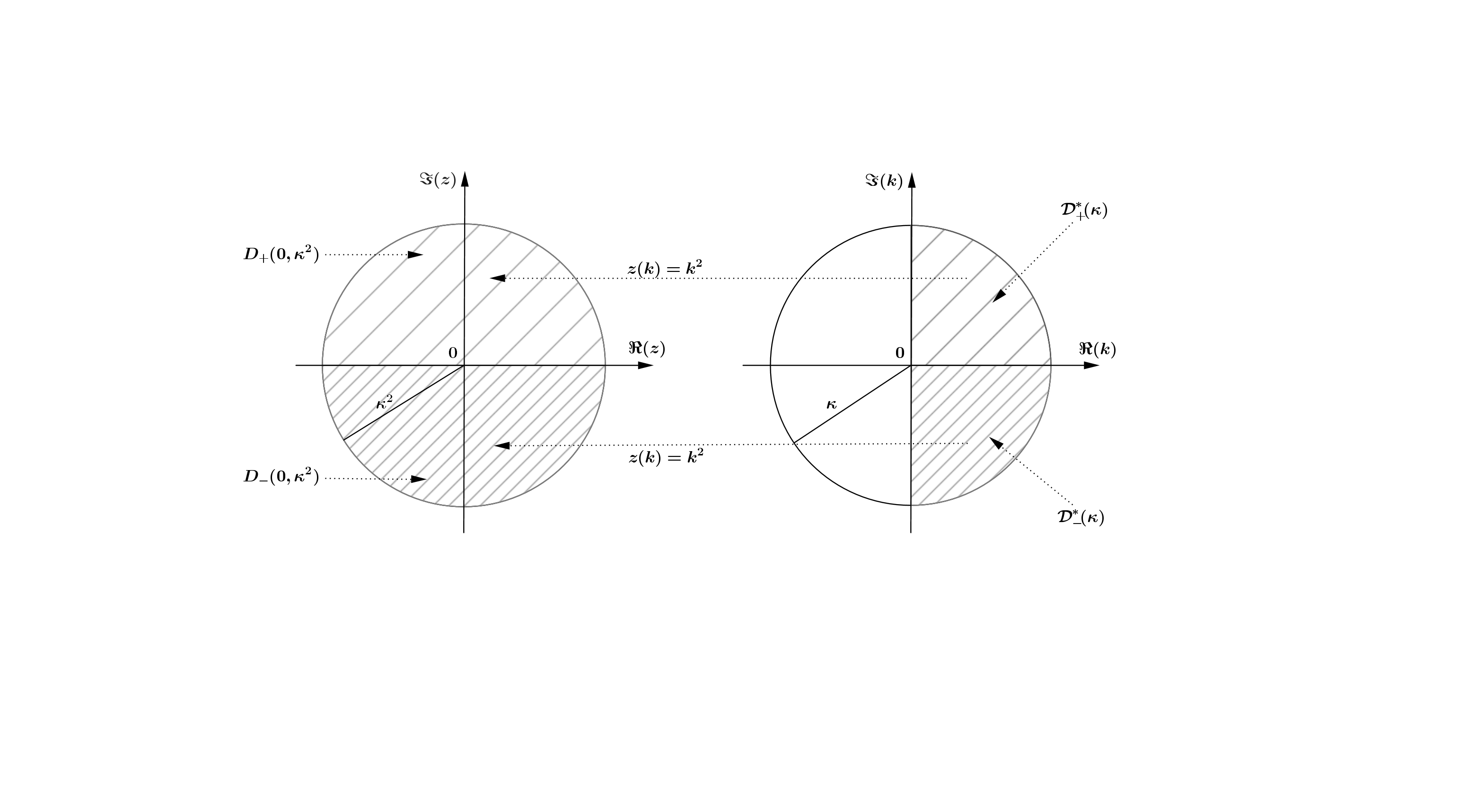}

\vspace*{-2cm}

\caption{Images $D_\pm(0,\kappa^2)$ of $\mathcal{D}_\pm^\ast(\kappa)$ 
by the local parametrisation $z(k) = k^2$.}
\end{center}
\end{figure}

Now, Assumption (C1) given by \eqref{eq1,13} implies 
that there exists a bounded operator $\mathscr{V}$ such that 
$\vert V \vert^\frac{1}{2} = \mathscr{V} G_\perp^\frac{1}{2}
G^\frac{1}{2}$, with $G_\perp^\frac{1}{2}$ and $G^\frac{1}{2}$ 
satisfying \eqref{eq1,13}. Then, the boundedness of $V$, identities \eqref{eq7,2}, 
\eqref{eq7,3}, together with Lemmas \ref{l7,1} and \ref{l7,2} give the following:

\begin{lem}\label{l7,3} 
The operator-valued functions
$$
\mathcal{D}_\pm^\ast (\kappa) \ni k \longmapsto 
\mathcal{T}_{V} \big( z(k) \big) := \Tilde{J} 
\vert V \vert^{\frac{1}{2}} \big( \htr(b,0) - z(k) \big)^{-1} 
\vert V \vert^{\frac{1}{2}}, \qquad z(k) := k^2,
$$ 
are holomorphic with values in $\sqq \big( L^{2}(\br^3) \big)$,
where $\Tilde{J}$ is defined by the polar decomposition
$V = \Tilde{J} \vert V \vert$ of $V$.
\end{lem}

\subsection{Reduction of the problem}

Similarly to Lemma \ref{l7,3}, we can show that
$V \big( \htr(b,0) - z(\cdot) \big)^{-1}$ is 
holomorphic in $\mathcal{D}_\pm^\ast(\kappa)$ with 
values in $\sqq \big( L^{2}(\br^3) \big)$. Then, as 
in \eqref{eq4,15}, we have
\begin{equation}\label{eq7,12}
z(k) \in \sigma_{\textup{\textbf{disc}}} \big( \htr(b,V) \big) 
\Leftrightarrow \textup{det}_{\lceil q \rceil} 
\left( I + V \big( \htr(b,0) - z(k) \big)^{-1} \right) = 0.
\end{equation}

\noindent
We are thus led to the following proposition:

\begin{prop}\label{p7,1} 
Let $\mathcal{T}_{V} \big( z(k) \big)$ be the operator 
defined in Lemma \ref{l7,3}. Then, the following 
assertions are equivalent:

\begin{itemize}
\item[(i)] $z(k) := k^2 \in D_\pm(0,\kappa^2)$ 
is a discrete eigenvalue of $\htr(b,V)$,

\item[(ii)] 
$\textup{det}_{\lceil q \rceil} \left( I + \mathcal{T}_{V} 
\big( z(k) \right) \big) = 0$,

\item[(iii)] $-1$ is an eigenvalue of 
$\mathcal{T}_{V} \big( z(k) \big)$.

\medskip

\noindent
Moreover,
\begin{equation}\label{eq7,13}
\textup{mult} \big( z(k) \big) = 
Ind_{\mathscr{C}} \hspace{0.5mm} 
\Big( I + \mathcal{T}_{V}\big( z(\cdot) \big) \Big),
\end{equation}
where $\mathscr{C}$ is a small contour positively 
oriented containing $k$ as the unique point 
$k \in \mathcal{D}_\pm^\ast(\kappa)$ verifying 
$z(k) \in D_\pm(0,\kappa^2)$ is a discrete 
eigenvalue of $\htr(b,V)$.
\end{itemize}
\end{prop}

\noindent
\begin{prof}
The proof is similar to that of Proposition 
\ref{p4,1}, taking into account the appropriate modifications.
\end{prof}

\subsection{Decomposition of the weighted resolvent}

Our goal in this section is to decompose $\mathcal{T}_{V} 
\big( z(k) \big)$ (defined just above), as a sum of a singular 
part at $k = 0$, and a holomorphic part in $\mathcal{D}_\pm^\ast(\kappa)$ 
continuous near $k = 0$ with values in $\sqq \big( L^{2}(\br^3) 
\big)$. The potential $V$ is supposed to verify Assumption 
(C1).

Observe that according to the choice \eqref{eq2,16} of the complex 
square root, we respectively have $\sqrt{k^2} = \pm k$ for $k \in 
\mathcal{D}_\pm^\ast(\kappa)$. By identity \eqref{eq7,3}, 
we have
\begin{equation}\label{eq7,14}
\mathcal{T}_{V} \big( z(k) \big) = \Tilde{J} 
\vert V \vert^{\frac{1}{2}} 
p \otimes \mathscr{R} \big( z(k) \big)
\begin{pmatrix}
1 & 0 \\
0 & 0
\end{pmatrix} \vert V \vert^{\frac{1}{2}}
+ \Tilde{J} \vert V \vert^{\frac{1}{2}} 
\big( \htr(b,0) - z(k) \big)^{-1} \textup{Q} 
\vert V \vert^{\frac{1}{2}}.
\end{equation}
Let us focus on the first term of the r.h.s. of 
\eqref{eq7,14} and define respectively $G_\pm$ as the 
multiplication operators by the functions 
$\br \ni X_\parallel \mapsto G^{\pm \frac{1}{2}}(X_\parallel)$. 
Then, we get
\begin{equation}\label{eq7,15}
\Tilde{J} 
\vert V \vert^{\frac{1}{2}} 
p \otimes \mathscr{R} \big( z(k) \big)
\begin{pmatrix}
1 & 0 \\
0 & 0
\end{pmatrix} \vert V \vert^{\frac{1}{2}}
= \Tilde{J} 
\vert V \vert^{\frac{1}{2}} G_-
p \otimes G_+ \mathscr{R} \big( z(k) \big) G_+
\begin{pmatrix}
1 & 0 \\
0 & 0
\end{pmatrix} G_- \vert V \vert^{\frac{1}{2}}.
\end{equation}
For $z \in \bc \setminus [0,+\infty)$,
$\mathscr{R}(z) := \big( -\frac{d^2}{dX_\parallel^2} - z)^{-1}$
admits the integral kernel
\begin{equation}\label{eq7,151}
I_z(X_\parallel,X_\parallel') := -\frac{e^{i\sqrt{z} 
\vert X_\parallel - X_\parallel' \vert}}{2i\sqrt{z}}.
\end{equation}  
Then, the integral kernel of the operator
$G_+ \mathscr{R} \big( z(k) \big) G_+$ is given by 
\begin{equation}\label{eq7,16}
\pm G^{\frac{1}{2}}(X_\parallel) 
\frac{i \textup{e}^{\pm i k \vert X_\parallel - X_\parallel' \vert}}{2 k}
G^{\frac{1}{2}}(X_\parallel'), \quad k \in \mathcal{D}_\pm^\ast(\kappa).
\end{equation}
With the help of \eqref{eq7,16}, we can write
\begin{equation}\label{eq7,17}
G_+ \mathscr{R} \big( z(k) \big) G_+ = 
\pm \frac{1}{k}a + b(k), \quad k \in 
\mathcal{D}_\pm^\ast(\kappa),
\end{equation}
$a : L^{2}(\mathbb{R}) \longrightarrow L^{2}(\mathbb{R})$ 
being the rank-one operator given by 
\begin{equation}\label{eq7,18}
a(u) := \frac{i}{2} \big\langle u,G^{\frac{1}{2}}(\cdot) 
\big\rangle G^{\frac{1}{2}}(X_\parallel),
\end{equation}
and $b(k)$ being the operator with integral kernel given by
\begin{equation}\label{eq7,19}
\pm G^{\frac{1}{2}}(X_\parallel) i \frac{ \textup{e}^{\pm i k 
\vert X_\parallel - X_\parallel' \vert} - 1}{2 k} 
G^{\frac{1}{2}}(X_\parallel').
\end{equation}
Moreover, it is easy to observe that we have $-2ia = c^\ast c$ 
with $c : L^{2}(\mathbb{R}) \longrightarrow \mathbb{C}$ 
defined by $c(u) := \langle u,G^\frac{1}{2}(\cdot) \rangle$, so 
that $c^{\ast} : \mathbb{C} \longrightarrow L^{2}(\mathbb{R})$ 
is given by $c^{\ast}(\lambda) = \lambda G^\frac{1}{2}(\cdot)$. 
Putting this together with \eqref{eq7,17} and \eqref{eq7,19},
we get
\begin{equation}\label{eq7,20}
\begin{aligned}
p \otimes G_+ \mathscr{R} \big( z(k) \big) G_+ 
= \pm \frac{i}{2k} p \otimes c^\ast c 
+ p \otimes s(k), \quad k \in \mathcal{D}_\pm^\ast(\kappa),
\end{aligned}
\end{equation}
where $s(k)$ is the operator acting from 
$G^{\frac{1}{2}}(X_\parallel) L^{2}(\mathbb{R})$ 
to $G^{-\frac{1}{2}}(X_\parallel) L^{2}(\mathbb{R})$ with integral 
kernel
\begin{equation}\label{eq7,21}
\pm \frac{ 1 - \textup{e}^{\pm i k \vert X_\parallel - X_\parallel' 
\vert}}{2 i k}.
\end{equation}
Equality \eqref{eq7,15} combined with \eqref{eq7,20} 
give for $k \in \mathcal{D}_\pm^\ast(\kappa)$
\begin{equation}\label{eq7,22}
\begin{split}
& \Tilde{J} 
\vert V \vert^{\frac{1}{2}} 
p \otimes \mathscr{R} \big( z(k) \big)
\begin{pmatrix}
1 & 0 \\
0 & 0
\end{pmatrix} \vert V \vert^{\frac{1}{2}} \\
& = \pm \frac{i\Tilde{J}}{2k} \vert V \vert^{\frac{1}{2}} 
G_- (p \otimes c^\ast c) \begin{pmatrix}
1 & 0 \\
0 & 0
\end{pmatrix} G_- \vert V \vert^{\frac{1}{2}} 
+ \Tilde{J} \vert V \vert^{\frac{1}{2}} G_- p \otimes 
s(k) 
\begin{pmatrix}
1 & 0 \\
0 & 0
\end{pmatrix}  G_- 
\vert V \vert^{\frac{1}{2}}.
\end{split}
\end{equation}
Finally, we get for $k \in \mathcal{D}_\pm^\ast(\kappa)$
\begin{equation}\label{eq7,23}
\Tilde{J} 
\vert V \vert^{\frac{1}{2}} 
p \otimes \mathscr{R} \big( z(k) \big)
\begin{pmatrix}
1 & 0 \\
0 & 0
\end{pmatrix} \vert V \vert^{\frac{1}{2}} = 
\pm \frac{i\Tilde{J}}{k} K^\ast K + \Tilde{J} 
\vert V \vert^{\frac{1}{2}} G_- p \otimes s(k) 
\begin{pmatrix}
1 & 0 \\
0 & 0
\end{pmatrix}
G_- \vert V \vert^{\frac{1}{2}},
\end{equation}
$K$ being defined by
\begin{equation}\label{eq7,24}
K := \frac{1}{\sqrt{2}} (p \otimes c) 
\begin{pmatrix}
1 & 0 \\
0 & 0
\end{pmatrix}
G_{-} \vert V \vert^{\frac{1}{2}}.
\end{equation} 
More precisely, recalling that $X_{\perp} := (x,y) \in \br^2$, 
we have $K : L^{2}(\br^{3}) \longrightarrow L^{2}(\br^{2})$ 
with 
$$
(K \psi)(X_{\perp}) = \frac{1}{\sqrt{2}} 
\int_{\br^{3}} dX_{\perp}^\prime dX_\parallel^\prime
{\mathcal P}_{b}(X_{\perp},X_{\perp}^\prime) 
\begin{pmatrix}
1 & 0 \\
0 & 0
\end{pmatrix}
\vert V \vert^{\frac{1}{2}} (X_{\perp}^\prime,X_\parallel^\prime) 
\psi (X_{\perp}^\prime,X_\parallel^\prime),
$$
where ${\mathcal P}_{b}(\cdot,\cdot)$ is the integral 
kernel of the orthogonal projection $p := p(b)$ (see \cite[Theorem 2.3]{hal}). 
Obviously, $K^{\ast} : L^{2}(\br^{2}) \longrightarrow L^{2}(\br^{3})$
is given by
$$
(K^{\ast}\varphi)(X_{\perp},X_\parallel) = 
\frac{1}{\sqrt{2}} \vert V \vert^{\frac{1}{2}} (X_\perp,X_\parallel) 
\begin{pmatrix}
1 & 0 \\
0 & 0
\end{pmatrix}
(p \varphi)(X_{\perp}).
$$
Therefore, it can be checked that the operator
$K K^{\ast} : L^{2}(\br^{2}) \longrightarrow L^{2}(\br^{2})$ 
satisfies
\begin{equation}\label{eq7,25}
K K^{\ast} = \begin{pmatrix}
1 & 0 \\
0 & 0
\end{pmatrix} 
p \textbf{\textup{V}}_{11} p,
\end{equation}
where $\textbf{\text{V}}_{11}$ is the multiplication 
operator by the function (also noted) $\textbf{\text{V}}_{11}$ 
defined by \eqref{eq2,14}.

For $\lambda \in \br_+ \setminus \lbrace 0 \rbrace$, we 
define $\big( -\frac{d^2}{dX_\parallel^2} - \lambda)^{-1}$ 
as the operator with integral kernel
\begin{equation}\label{eq7,251}
\displaystyle I_\lambda(X_\parallel,X_\parallel') := 
\lim_{\delta \downarrow 0} I_{\lambda + i\delta} 
(X_\parallel,X_\parallel') =
\frac{ie^{i\sqrt{\lambda}\vert X_\parallel - X_\parallel' 
\vert}}{2\sqrt{\lambda}},
\end{equation}
$I_z(\cdot)$ being defined by \eqref{eq7,151}. Then, as in
\cite[(Proof of) Proposition 4.2]{rage}, we can show by a limiting 
absorption principle that the operator-valued function 
$\overline{\mathcal{D}_\pm^\ast(\kappa)} \ni k \mapsto 
G_+ s(k) G_+ \in \sd \big( L^{2}(\br) \big)$ 
is well defined and continuous. We thus have proved the
following:

\begin{prop}\label{p7,2} 
For $k \in \mathcal{D}_\pm^\ast(\kappa)$, we have
\begin{equation}\label{eq7,26}
\mathcal{T}_{V} \big( z(k) \big) = 
\pm \frac{i\Tilde{J}}{k} \mathscr{B}
+ \mathscr{A}(k), \quad \mathscr{B} := K^\ast K,
\end{equation}
where the operator $\mathscr{A}(k) \in 
\sqq \big( L^{2}(\br^3) \big)$ given by
$$
\mathscr{A}(k) := 
\Tilde{J} 
\vert V \vert^{\frac{1}{2}} G_- p \otimes s(k) 
\begin{pmatrix}
1 & 0 \\
0 & 0
\end{pmatrix}
G_- \vert V \vert^{\frac{1}{2}}
 + \Tilde{J} \vert V \vert^{\frac{1}{2}} 
\big( \htr(b,0) - z(k) \big)^{-1} \textup{Q} 
\vert V \vert^{\frac{1}{2}},
$$
is holomorphic in $\mathcal{D}_\pm^\ast(\kappa)$ and
continuous on $\overline{\mathcal{D}_\pm^\ast(\kappa)}$,
with $s(k)$ defined by \eqref{eq7,20}.
\end{prop} 

\begin{rem}\label{r7,1}
\end{rem}

\begin{itemize}
\item[(i)] For any $r > 0$, according to \eqref{eq7,25}, 
we have
\begin{equation}\label{eq7,27}
\textup{Tr} \hspace{0.4mm} \one_{(r,\infty)} 
\left( K^\ast K \right) 
= \textup{Tr} \hspace{0.4mm} \one_{(r,\infty)} 
\left( K K^\ast \right)
= \textup{Tr} \hspace{0.4mm} \one_{(r,\infty)} 
\big( p \textbf{\textup{V}}_{11} p \big).
\end{equation}

\item[(ii)] If $V$ satisfies Assumption (C2) given by 
\eqref{eq2,20}, then Proposition \ref{p4,2} holds with 
$\Tilde{J}$ replaced by $J\eta$, $J := sign(W)$, and 
$\textbf{\textup{V}}_{11}$ replaced by 
$\textbf{\textup{W}}_{11}$ in \eqref{eq7,27}.
\end{itemize}

\section{Proof of Theorem \ref{t2,4}: Upper bounds}\label{s8}

The proof is similar to that of Theorem \ref{t2,1}.

\subsection{A preliminary result}

Introduce the numerical range
$$
N \big( \htr(b,V) \big) := \big\lbrace \langle 
\htr(b,V)f,f \rangle : f \in Dom \big( \htr(b,V) \big), 
\Vert f \Vert_{L^{2}} = 1 \big\rbrace,
$$ 
satisfying $\sigma \big( \htr(b,V) \big) \subseteq 
\overline{N \big( \htr(b,V) \big)}$.

\begin{prop}\label{p8,1} 
There exists $r_{0} > 0$ such that for any $k \in 
\big\lbrace 0 < r < \vert k \vert < r_{0} \big\rbrace \cap 
\mathcal{D}_\pm^\ast(\kappa)$, we have:
 
\begin{itemize}
\item[(i)] $z(k) := k^2$ is a discrete eigenvalue of 
$\htr(b,V)$ near zero if and only if $k$ is a zero of
\begin{equation}\label{eq8,1}
\mathscr{D}(k,r) := \det \big( I + \mathscr{K}(k,r) \big),
\end{equation}
with $\mathscr{K}(k,r)$ a finite-rank operator analytic 
with respect to $k$ and satisfying
$$
\small{\textup{rank} \hspace{0.6mm} \mathscr{K}(k,r) = 
\mathcal{O} \Big( \textup{Tr} \hspace{0.4mm} 
\one_{(r,\infty)} 
\big( p \textbf{\textup{V}}_{11} p \big) + 1 \Big)},
\quad \left\Vert \mathscr{K}(k,r) \right\Vert = 
\mathcal{O} \left( r^{-1} \right),
$$
where the $\mathcal{O}$'s are uniform with respect to $r$, $k$.

\item[(ii)] Furthermore, if $z(k) := k^2$ is a discrete 
eigenvalue of $\htr(b,V)$ near zero, then we have
\begin{equation}\label{eq8,2}
\textup{mult} \big( z(k) \big) = Ind_{\mathscr{C}} 
\hspace{0.5mm} \left( I + \mathscr{K}(\cdot,r) \right) 
= \textup{m}(k),
\end{equation}
$\mathscr{C}$ being chosen as in \eqref{eq7,13}, and 
$\textup{m}(k)$ being the multiplicity of $k$ as a
zero of $\mathscr{D}(\cdot,r)$.

\item[(iii)] If $z(k)$ satisfies 
$\textup{dist} \big( z(k),\overline{N \big( \htr(b,V) \big)} 
\big) > \varsigma > 0$, $\varsigma = \mathcal{O}(1)$, then
$I + \mathscr{K}(k,r)$ is invertible and satisfies
$
\left\Vert \big( I + \mathscr{K}(k,r) \big)^{-1} \right\Vert = 
\mathcal{O} \left( \varsigma^{-1} \right)
$,
where the $\mathcal{O}$ is uniform with respect to $r$, $k$ and
$\varsigma$.
\end{itemize}
\end{prop}

\noindent
\begin{prof}
The proof follows by arguing similarly to that of Proposition 
\ref{p5,1}, taking into account the appropriate modifications.
\end{prof}

\subsection{Back tot the proof Theorem \ref{t2,4}}

Proposition \ref{p8,1} above implies that
\begin{equation}\label{eq8,3}
\begin{aligned}
\mathscr{D}(k,r) & = \prod_{j=1}^{\mathcal{O} \big( \textup{Tr} 
\hspace{0.4mm} \one_{(r,\infty)} 
(p \textbf{\textup{V}}_{11} p) + 1 \big)} 
\big{(} 1 + \lambda_{j}(k,r) \big{)}\\
& = \mathcal{O}(1) \hspace{0.5mm} \textup{exp} \hspace{0.5mm} 
\Big( \mathcal{O} \big( \textup{Tr} \hspace{0.4mm} 
\one_{(r,\infty)} 
\big( p \textbf{\textup{V}}_{11} p \big) + 1 \big) 
\vert \ln r \vert \Big),
\end{aligned}
\end{equation}
for $0 < r < \vert k \vert < r_{0}$, the $\lambda_{j}(k,r)$ 
being the eigenvalues of $\mathscr{K} := \mathscr{K}(k,r)$ 
satisfying $\vert \lambda_{j}(k,r) \vert = \mathcal{O} 
\left( r^{-1} \right)$. If
$\textup{dist} \big( z(k),\overline{N \big( \htr(b,V) \big)} 
\big) > \varsigma > 0$ with $0 < r < \vert k \vert < r_{0}$,
then,
$$
\mathscr{D}(k,r)^{-1} = \det \big( I + \mathscr{K} \big)^{-1} 
= \det \big( I - \mathscr{K} ( I + \mathscr{K})^{-1} \big).
$$
Similarly to \eqref{eq8,3}, we can show that
\begin{equation}\label{eq8,4}
\vert \mathscr{D}(k,r) \vert \geq C \hspace{0.5mm} 
\textup{exp} \hspace{0.5mm} 
\Big( - C \big( \textup{Tr} \hspace{0.4mm} 
\one_{(r,\infty)} 
\big( p \textbf{\textup{V}}_{11} p \big) + 1 \big) 
\big( \vert \textup{ln} \hspace{0.5mm} \varsigma \vert + 
\vert \textup{ln} \hspace{0.5mm} r \vert \big) \Big),
\end{equation}
so that for $r^2 < \varsigma < 4r^2$, $0 < r \ll 1$, we obtain
\begin{equation}\label{eq8,40}
- \ln \, \vert \mathscr{D}(k,r) \vert \leq
C \, \textup{Tr} \hspace{0.4mm} \one_{(r,\infty)} 
\big( p \vert V \vert_{11} p \big) \vert \textup{ln} \hspace{0.5mm} r \vert
+ \mathcal{O}(1).
\end{equation}
Consider the domains $\Delta_\pm := \left\lbrace r < \vert k 
\vert < 2r : \vert \Re(k) \vert > \sqrt{\frac{\nu}{2}} : 
\vert \Im(k) \vert > \sqrt{\frac{\nu}{2}} \right\rbrace \cap 
\mathcal{D}_\pm^\ast(\kappa)$ with 
$0 < r < \sqrt{\Vert V \Vert} < \sqrt{\frac{5}{2}}r$ and 
$0 < \nu < 2r^2$. Since the numerical range of the operator $\htr(b,V)$
satisfies
\begin{equation}
N \big( \htr(b,V) \big) \subseteq \big\lbrace z \in \bc : \vert \Im(z) 
\vert \leq \Vert V \Vert \big\rbrace,
\end{equation}
then there exists some $k_{0} \in \Delta_\pm /r$ satisfying 
$\textup{dist} \big( z(rk_0),\overline{N \big( \htr(b,V) \big)} 
\big) \geq \varsigma > r^2$, $\varsigma < 4r^2$. Then, Theorem \ref{t2,4} 
follows by applying the Jensen Lemma \ref{la,1} with the function 
$g(k) := \mathscr{D}(rk,r)$, $k \in \Delta_\pm /r$, together with 
\eqref{eq8,3} and \eqref{eq8,40}.

\section{Proof of Theorem \ref{t2,6}: Sectors free of complex eigenvalues and lower bounds}\label{s10}

To simplify, we give the proof only for the case 
$\alpha \in (0,\pi)$. The case $\alpha \in 
-(0,\pi)$ follows in a similar way by replacing $k$ by $-k$.

\medskip

\textit{(i):} For any $\theta > 0$ small enough, set $\delta = \tan(\theta)$
and introduce the sector 
\begin{equation}\label{eq2,21}
\mathcal{C}_\delta := \big\lbrace k \in \bc : 
- \delta \Im(k) \leq \vert \Re(k) \vert \big\rbrace.
\end{equation}
According to (ii) of Remark \ref{r7,1}, for any 
$\varepsilon > 0$, we have
\begin{equation}\label{eq10,1}
I + \mathcal{T}_{\varepsilon V} \big( z(k) \big) = 
I + \frac{i \varepsilon \eta}{k} \mathscr{B} + 
\varepsilon \mathscr{A}(k), \qquad
k \in \mathcal{D}_+^\ast(\kappa),
\end{equation}
$\mathscr{B}$ being a self-adjoint positive operator 
independent of $k$, while $\mathscr{A}(k) \in \sqq \big( 
L^{2}(\br^3) \big)$ is holomorphic in 
$\mathcal{D}_+^\ast(\kappa)$ and continuous on 
$\overline{\mathcal{D}_+^\ast(\kappa)}$. Since we have
$$
I + \frac{i \varepsilon \eta}{k} \mathscr{B} = 
\frac{i \eta}{k} (\varepsilon \mathscr{B} - ik \eta^{-1}),
$$
then it is easy to see that the operator $I + \frac{i 
\varepsilon \eta}{k} \mathscr{B}$ is invertible for 
$ik \eta^{-1} \notin \sigma (\varepsilon \mathscr{B})$. 
Otherwise, it can be shown that we have
\begin{equation}\label{eq9,2}
\left\Vert \left( I + \frac{i \varepsilon \eta}{k} 
\mathscr{B} \right)^{-1} \right\Vert \leq 
\frac{\vert k \eta^{-1} \vert}{\sqrt{\big( \Im(k \eta^{-1}) 
\big)_+^2 + \vert \Re(k \eta^{-1}) \vert^2}}, 
\qquad r_+ := \max (r,0).
\end{equation}
Therefore, for $k \in \eta \mathcal{C}_\delta$, 
it can be checked that 
\begin{equation}\label{eq9,3}
\left\Vert \left( I + \frac{i \varepsilon \eta}{k} 
\mathscr{B} \right)^{-1} \right\Vert \leq \sqrt{1 + \delta^{-2}},
\end{equation}
$k \in \mathcal{D}_+^\ast(\kappa)$. Then, we have 
\begin{equation}\label{eq10,3}
I + \mathcal{T}_{\varepsilon V} \big( z(k) \big) = 
\big( I + A(k) \big) \left( I + \frac{i\varepsilon 
\eta}{k} \mathscr{B} \right),
\end{equation}
where
\begin{equation}\label{eq10,4}
A(k) := \varepsilon \mathscr{A} (k)
\left( I + \frac{i\varepsilon \eta}{k} \mathscr{B} 
\right)^{-1} \in \sqq \big( L^{2}(\br^3) \big).
\end{equation}
Since $\mathscr{A}(k) \in \sqq \big( L^{2}(\br^3) \big)$ is 
continuous on $\overline{\mathcal{D}_+^\ast(\kappa)}$, then 
there exists a uniform constant $C > 0$ such that 
$\Vert \mathscr{A}(k) \Vert \le \Vert \mathscr{A}(k) \Vert_{\sqq} \leq C$. 
Putting this together with \eqref{eq9,3} and \eqref{eq10,4}, it 
follows immediately that 
$I + \mathcal{T}_{\varepsilon V} \big( z(k) \big)$ 
is invertible for $k \in \eta \mathcal{C}_\delta$,
$k \in \mathcal{D}_+^\ast(\kappa)$,
and 
\begin{equation}\label{eq10,40}
0 < \varepsilon < C_0 := \big( C \sqrt{1 + \delta^{-2}} \big)^{-1}.
\end{equation} 
This means that $z(k)$ is not a discrete eigenvalue.

\medskip

\textit{(ii):} Let $(\mu_j)_j$ denote the sequence of the decreasing 
nonzero eigenvalues of $p\textbf{\textup{W}}_{11}p$ 
taking into account their multiplicity. If Assumption (C3) 
given by \eqref{eq2,25} is fulfilled, then similarly to 
\cite[(Proof of) Lemma 7]{bon}, it can be shown 
that there exists a positive constant $\nu$ such that
\begin{equation}\label{eq10,16}
\# \big\lbrace j : \mu_j - \mu_{j+1} > \nu \mu_j \big\rbrace 
= \infty.
\end{equation}
Since the nonzero eigenvalues of $\mathscr{B}$ and 
$p\textbf{\textup{W}}_{11}p$ coincide, then
there exists a decreasing sequence of positive numbers 
$(r_\ell)_\ell$, $r_{\ell} \searrow 0$, such that
\begin{equation}\label{eq10,17}
\textup{dist} \big( r_\ell,\sigma(\mathscr{B}) \big) 
\geq \frac{\nu r_\ell}{2}, \quad \ell \in \bn.
\end{equation}
Moreover, there exists for any $\ell \in \bn$ a path 
$\Tilde{\Sigma}_{\ell} := \partial 
\Lambda_\ell$, where
\begin{equation}\label{eq10,181}
\Lambda_\ell := \big\lbrace \Tilde{k} \in \bc :
0 < \vert \Tilde{k} \vert < r_0 : \vert \Im(\Tilde{k}) \vert \leq 
\delta \Re(\Tilde{k}) : r_{\ell + 1} \leq \Re(\Tilde{k}) \leq r_\ell 
\big\rbrace,
\end{equation}
(see Figure 8.1), enclosing the eigenvalues of $\mathscr{B}$ lying in 
$[r_{\ell + 1},r_\ell]$.

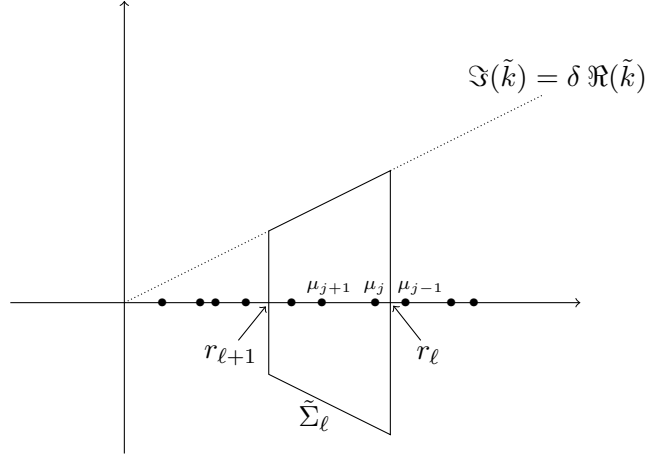
\begin{figure}[h]\label{fig 4}
\begin{center}
\tikzstyle{+grisEncadre}=[dashed]
\tikzstyle{blancEncadre}=[fill=white!100]
\tikzstyle{grisEncadre}=[densely dotted]
\tikzstyle{dEncadre}=[dotted]

\begin{tikzpicture}[scale=1]

\draw [->] (0,-2) -- (0,4);
\draw [->] (-1.5,0) -- (6,0);

\draw (1.9,-0.95) -- (1.9,0.95) -- (3.5,1.75) -- (3.5,-1.75) -- cycle;

\draw [grisEncadre] (0,0) -- (1.9,0.95);
\draw [grisEncadre] (3.5,1.75) -- (5.5,2.75);

\node at (5.7,3) {$\Im(\Tilde{k}) = \delta \hspace*{0.08cm} \Re(\Tilde{k})$};
\node at (2.5,-1.5) {$\Tilde{\Sigma}_{\ell}$};

\draw [->] (1.5,-0.5) -- (1.86,-0.05);
\node at (1.4,-0.7) {$r_{\ell + 1}$};

\draw [->] (3.9,-0.5) -- (3.54,-0.04);
\node at (4,-0.7) {$r_{\ell}$};

\node at (3.3,0) {\tiny{$\bullet$}};
\node at (3.3,0.2) {\tiny{$\mu_{j}$}};

\node at (3.7,0) {\tiny{$\bullet$}};
\node at (3.9,0.2) {\tiny{$\mu_{j-1}$}};

\node at (1.6,0) {\tiny{$\bullet$}};
\node at (2.2,0) {\tiny{$\bullet$}};

\node at (2.6,0) {\tiny{$\bullet$}};
\node at (2.7,0.2) {\tiny{$\mu_{j+1}$}};

\node at (1.2,0) {\tiny{$\bullet$}};
\node at (1,0) {\tiny{$\bullet$}};
\node at (0.5,0) {\tiny{$\bullet$}};

\node at (4.3,0) {\tiny{$\bullet$}};
\node at (4.6,0) {\tiny{$\bullet$}};

\end{tikzpicture}
\caption{Representation of the path $\Tilde{\Sigma}_\ell  
= \partial \Lambda_\ell$.}
\end{center}
\end{figure}

\noindent
Obviously, the operator $\Tilde{k} - \mathscr{B}$ is invertible 
for $\Tilde{k} \in \Tilde{\Sigma}_{\ell}$, and it can be proved that 
\begin{equation}\label{eq10,18}
\big\Vert (\Tilde{k} - \mathscr{B})^{-1} \big\Vert 
\le \frac{1}{{\rm dist} \big( \Tilde{k},\sigma(\mathscr{B}) \big)}
= \frac{1}{ \vert \Tilde{k} \vert} \times \frac{\vert \Tilde{k} \vert}{{\rm dist} \big( \Tilde{k},\sigma(\mathscr{B}) \big)}
\le \frac{C(\delta,\nu)}{ \vert \Tilde{k} \vert},
\end{equation}
where 
\begin{equation}\label{eq10,180}
C(\delta,\nu) := \sqrt{1 + \delta^2} \max \left( \delta^{-1},(\nu/2)^{-1} \right).
\end{equation} 
Introduce the path $\Sigma_{\ell} := -i\varepsilon \eta 
\Tilde{\Sigma}_{\ell}$.
According to the construction of $\Sigma_{\ell}$ and 
\eqref{eq10,18}, we immediately observe that $I + 
\frac{i\varepsilon \eta}{k} \mathscr{B}$ is invertible for 
$k \in \Sigma_{\ell}$ with
\begin{equation}\label{eq10,20}
\left\Vert \left( I + \frac{i\varepsilon \eta}{k} 
\mathscr{B} \right)^{-1} \right\Vert
\leq C(\delta,\nu).
\end{equation}
Then, for $k \in \Sigma_{\ell}$, we have
\begin{equation}\label{eq10,21}
I + \frac{i\varepsilon \eta}{k} \mathscr{B} + 
\varepsilon \mathscr{A}(k) = \left( I + \varepsilon 
\mathscr{A}(k) \left( I + \frac{i\varepsilon \eta}{k} 
\mathscr{B} \right)^{-1} \right) \left( I + 
\frac{i\varepsilon \eta}{k} \mathscr{B} \right).
\end{equation}
By choosing $0 < \varepsilon \leq \varepsilon_0$ sufficiently
small and using Property \textbf{e)} of Subsection \ref{s3,1} 
given by \eqref{eq3,5}, we obtain
\begin{equation}\label{eq7,22}
\left\vert \textup{det}_{\lceil q \rceil} \left[ I + \varepsilon 
\mathscr{A}(k) \left( I + \frac{i \varepsilon \eta}{k} 
\mathscr{B} \right)^{-1} \right] - 1 \right\vert < 1,
\end{equation}
for any $k \in \Sigma_{\ell}$. More precisely, if we let $C$, $C_0$ 
be the constants defined by \eqref{eq10,40}, $C(\delta,\nu)$ the one 
defined by \eqref{eq10,180}, and $\Gamma_q$ that defined by \eqref{eq3,5}, 
then \eqref{eq7,22} holds whenever $\varepsilon$ satisfies
\begin{equation}\label{eq7,220}
0 < \varepsilon < C^{-1} C(\delta,\nu)^{-1} e^{-\Gamma_q \big( C_0 C C(\delta,\nu)  
+ 1 \big)^{\lceil q \rceil}} = C_0 C_1(\delta,\nu)^{-1} e^{-\Gamma_q \big( C_1(\delta,\nu)  
+ 1 \big)^{\lceil q \rceil}},
\end{equation}
where
\begin{equation}\label{eq7,221}
C_1(\delta,\nu) := \delta \max \left( \delta^{-1},(\nu/2)^{-1} \right).
\end{equation} 
Thus, the Rouché Theorem implies that the number of zeros of 
$\textup{det}_{\lceil q \rceil} 
\big( I + \frac{i\varepsilon \eta}{k} \mathscr{B} 
+ \varepsilon \mathscr{A}(k) \big)$ enclosed in
$\big\lbrace z(k) \in D_+(0,\kappa^2) : k \in 
\Lambda_{\ell} \big\rbrace$ taking into account their 
multiplicity, is equal to that of 
$\textup{det}_{\lceil q \rceil} \big( I + \frac{i \varepsilon 
\eta}{k} \mathscr{B} \big)$ enclosed in $\big\lbrace z(k) \in 
D_+(0,\kappa^2) : k \in \Lambda_{\ell} \big\rbrace$ taking 
into account their multiplicity. This number is equal to
$\textup{Tr} \hspace{0.4mm} \one_{[r_{\ell +1},r_\ell]} \big( p 
\textbf{\textup{W}}_{11}p \big)$. Hence, thanks to
Proposition \ref{p7,1} and Property \eqref{eqa,3} applied 
to \eqref{eq10,21}, bound \eqref{lb1} follows immediately 
since the zeros of $\textup{det}_{\lceil q \rceil} \big( I + 
\frac{i\varepsilon \eta}{k} \mathscr{B} + \varepsilon 
\mathscr{A}(k) \big)$ are the discrete eigenvalues of 
$\htr(b,\varepsilon V)$ taking into account their 
multiplicity. From the fact that the sequence $(r_\ell)_\ell$ 
is infinite tending to zero, it follows the infiniteness 
of the number of the discrete eigenvalues claimed. 
This concludes the proof of Theorem $\ref{t2,6}$.

\section{Appendix A $1$: Schatten-von Neumann ideals and regularized determinants}\label{s3,1}

For the convenience of the reader, we repeat the relevant 
material from Reed-Simon \cite{ree}, Simon \cite{simo,sim}, and 
Gohberg-Goldberg-Krupnik \cite{gohb}, thus making our exposition
self-contained.

Let $\mathscr{H}$ be a separable Hilbert space and $\sinf(\mathscr{H})$ 
be the set of compact linear operators on $\mathscr{H}$. Denote 
by $s_k(T)$ the $k$-th singular value of $T \in \sinf(\mathscr{H})$. 
The Schatten-von Neumann classes
are defined by 
\begin{equation}\label{eq3,1}
\sqq(\mathscr{H}) := \Big\lbrace T \in \sinf(\mathscr{H}) : 
\Vert T \Vert^q_\sqq := \sum_k s_k(T)^q < +\infty \Big\rbrace,
\quad q \in [1,+\infty).
\end{equation}
To simplify, we will write $\sqq$ when no confusion can arise.
For  $\lceil q \rceil := \min \big\lbrace n \in \mathbb{N} : n \geq q 
\big\rbrace$ and $T \in \sqq$, the regularized determinant is defined by
\begin{equation}\label{eq3,2}
\textup{det}_{\lceil q \rceil} (I - T)
 := \prod_{\mu \hspace*{0.1cm} \in \hspace*{0.1cm} \sigma (T)} 
 \left[ (1 - \mu) \exp \left( \sum_{k=1}^{\lceil q \rceil-1} \frac{\mu^{k}}{k} 
 \right) \right].
\end{equation}
Here are some elementary properties about this determinant 
(see for instance \cite{simo}):

\textbf{a)} $\textup{det}_{\lceil q \rceil} (I) = 1$.

\textbf{b)} For $A$, $B \in \mathscr{L} (\mathscr{H})$ the class of bounded
linear operators on $\mathscr{H}$, if $AB$ and $BA$ belong to $\sqq$, then
$\textup{det}_{\lceil q \rceil} (I - AB)
= \textup{det}_{\lceil q \rceil} (I - BA)$.

\textbf{c)} $I - T$ is invertible if and only if
$\textup{det}_{\lceil q \rceil} (I - T) \neq 0$.

\textbf{d)} If $T : D \longrightarrow \sqq$ is a holomorphic
operator-valued function in a domain $D$, then so is
$\textup{det}_{\lceil q \rceil} \big( I - T(\cdot) \big)$ in $D$.

\textbf{e)} $\textup{det}_{\lceil q \rceil} (I - T)$ is Lipschitz as 
function on $\sqq$ uniformly on balls. Explicitly, we have
\begin{equation}\label{eq3,5}
\big\vert \textup{det}_{\lceil q \rceil} (I - T_1) - 
\textup{det}_{\lceil q \rceil} (I - T_2) \big\vert
\leq \Vert T_1 - T_2 \Vert_\sqq 
e^{\Gamma_q \big( \Vert T_1 \Vert_\sqq + \Vert T_2 \Vert_\sqq + 1 \big)^{\lceil q \rceil}},
\end{equation}
by \cite[Theorem 6.5]{simo}, for some constant $\Gamma_q > 0$.

\section{Appendix A $2$: Index of a finite meromorphic operator-valued function}\label{sa,1}


The space $\mathscr{H}$ and the class $\mathscr{L}(\mathscr{H})$ are defined as in 
Appendix A 1. We have the following definition from \cite[Definition 4.1.1]{goh}.

\begin{fe}
Let $\mathcal{U}$ be a neighbourhood of a fixed point $w \in \bc$, and 
$F : \mathcal{U} \setminus \lbrace w \rbrace \longrightarrow \mathscr{L}(\mathscr{H})$ 
be a holomorphic operator-valued function. The function $F$ is said to be finite 
meromorphic at $w$ if its Laurent expansion at $w$ has the form
\begin{equation}
F(z) = \sum_{n = m}^{+\infty} (z - w)^n A_n, \quad m > - \infty,
\end{equation}
where for $m < 0$, the operators $A_m, \ldots, A_{-1}$ are of finite rank.
Moreover, if $A_0$ is a Fredholm operator, then the function $F$ is said to be Fredholm 
at $w$. In that case, the Fredholm index of $A_0$ is called the Fredholm index of $F$ 
at $w$.
\end{fe}

If a function $f$ is holomorphic in a neighbourhood of a contour 
$\mathscr{C}$ (positively oriented), its index with respect to this
contour is defined by 
\begin{equation}\label{eqa,1}
ind_{\mathscr{C}} \hspace{0.5mm} f 
:= \frac{1}{2i\pi} \int_{\mathscr{C}} \frac{f'(z)}{f(z)} dz.
\end{equation}
Let us point out that if $f$ is holomorphic in a domain $D$ 
with $\partial D = \mathscr{C}$, then thanks to the residues 
theorem, $\textup{ind}_{\mathscr{C}} \hspace{0.5mm} f$ coincides 
with the number of zeros of $f$ in $D$ taking into account 
their multiplicity. 

In what follows below, $\textup{GL}(\mathscr{H})$ denotes the class of 
invertible linear operators on the Hilbert space $\mathscr{H}$. 
Let $D \subseteq \mathbb{C}$ be a connected domain, 
$Z \subset D$ be a pure point and closed 
subset, and $A : \overline{D} \backslash Z \longrightarrow 
\textup{GL}(\mathscr{H})$
be a finite meromorphic operator-valued function which is Fredholm 
at each point of $Z$. The index of $A$ with respect 
to the contour $\partial D$ is defined by 
\begin{equation}\label{eqa,2}
Ind_{\partial D} \hspace{0.5mm} A := 
\frac{1}{2i\pi} \textup{Tr} \int_{\partial D} A'(z)A(z)^{-1} 
dz = \frac{1}{2i\pi} \textup{Tr} \int_{\partial D} A(z)^{-1} 
A'(z) dz,
\end{equation} 
where the operator $A$ does not vanish in the integration contour 
$\partial D$. The following properties are well known: 
\begin{equation}\label{eqa,3}
Ind_{\partial D} \hspace{0.5mm} A_{1} A_{2} = 
Ind_{\partial D} \hspace{0.5mm} A_{1} + 
Ind_{\partial D} \hspace{0.5mm} A_{2};
\end{equation} 
for $K(z)$ a trace class operator-valued function, we have
\begin{equation}\label{eqa,4}
Ind_{\partial D} \hspace{0.5mm} (I+K)= 
ind_{\partial D} \hspace{0.5mm} \det \hspace{0.5mm} (I + K).
\end{equation} 
We refer for instance to \cite[Chap. 4]{goh} for a deeper discussion on
the subject.

\section{Appendix A $3$: Jensen type inequality and characteristic values of operator-valued functions}\label{sa,3}

The following lemma \big(see for instance \cite[Lemma 6]{bon} 
for a proof\big) contains a version of the well-known Jensen 
inequality.

\begin{lem}\label{la,1} 
Let $\Delta$ be a simply connected 
sub-domain of $\mathbb{C}$ and let $g$ be holomorphic 
in $\Delta$ with continuous extension 
to $\overline{\Delta}$. Assume that there exists 
$\lambda_{0} \in \Delta$ such that $g(\lambda_{0}) \neq 0$ 
and $g(\lambda) \neq 0$ for $\lambda\in \partial \Delta$
(the boundary of $\Delta$). 
Let $\lambda_{1}, \lambda_{2}, \ldots, \lambda_{N} \in \Delta$ 
be the zeros of $g$ repeated according to their multiplicity. 
For any domain $\Delta' \subset \subset \Delta$, there exists 
$C' > 0$ such that $N(\Delta',g)$, the number of zeros 
$\lambda_{j}$ of $g$ contained in $\Delta'$, satisfies
\begin{equation}\label{eqa,5}
N(\Delta',g) \leq C' \left( \int_{\partial \Delta} 
\textup{ln} \vert g(\lambda) \vert d\lambda 
- \textup{ln} \vert g(\lambda_{0}) \vert  \right).
\end{equation}
\end{lem}

Consider a domain $\mathcal{D}$ of $\mathbb{C}$ containing $0$,
and let $T : \mathcal{D} \longrightarrow S_{\infty}(\mathscr{H})$ 
be a holomorphic operator-valued function, $\mathscr{H}$ being as
above.

\begin{fe}\label{d1}
For a domain 
$\Omega \subset \mathcal{D} \setminus \lbrace 0 \rbrace$, a 
complex number $z \in \Omega$ is a characteristic value of 
$z \mapsto \mathscr{T}(z) := I - \frac{T(z)}{z}$, if
$\mathscr{T}(z)$ is not invertible. The multiplicity 
of a characteristic value $z_{0}$ is defined by
\begin{equation}\label{eqa,6}
\textup{mult}(z_{0}) := Ind_{\mathscr{C}} \big( I - \mathscr{T}(\cdot) 
\big),
\end{equation}
where $\mathscr{C}$ is a small contour positively oriented 
containing $z_{0}$ as the unique point $z$ satisfying 
$\mathscr{T}(z)$ is not invertible.
\end{fe}

\noindent
Define
$$
\mathcal{Z}(\Omega) := \big\lbrace z \in \Omega : 
\mathscr{T}(z) \hspace{0.8mm} \textup{is not invertible} 
\big\rbrace.
$$ 
Once there exists $z_{0} \in \Omega$ satisfying $\mathscr{T}(z_{0})$ 
is not invertible, then by the analytic Fredholm theorem, 
the set $\mathcal{Z}(\Omega)$ is pure point. Hence, we set
\begin{equation}\label{eqa,7}
\mathcal{N}(\Omega) := \# \mathcal{Z}(\Omega).
\end{equation}
In the sequel, we suppose that the operator $T(0)$ is 
self-adjoint and we put 
\begin{equation}\label{eqa,8}
n(\omega) := \textup{Tr} \hspace{0.6mm} \one_{\omega} 
\big( T(0) \big),
\end{equation}
the number of eigenvalues of $T(0)$ lying in the interval 
$\omega \subset \mathbb{R}^{\ast}$, taking into account their 
multiplicity. The orthogonal projection onto 
$\textup{Ker} \hspace{0.6mm} T(0)$ is denoted $\Pi_{0}$.

\begin{lem}\label{la,2} 
\cite[Corollary 3.4, Corollary 3.9, Corollary 3.11]{bo} 
Let $T$ be as above with $I - T'(0) \Pi_{0}$ invertible. Let 
$\Omega \subset \bc \setminus \lbrace 0 \rbrace$ be a bounded 
domain such that  $\partial \Omega$ is smooth and transverse 
to the real axis at each point of $\partial \Omega \cap \br$.

\begin{itemize}
\item[(i)] If $\Omega \cap \mathbb{R} = \emptyset$,
then for $s$ sufficiently small, $\mathcal{N}(s\Omega) = 0$. 
So, the characteristic values $z \in \mathcal{Z}(\Omega)$ 
satisfy $\vert \Im(z) \vert = o(\vert z \vert)$ near $0$.

\item[(ii)] Moreover, if $T(0)$ satisfy $\pm T(0) 
\geq 0$, then the characteristic values $z$ satisfy respectively 
$\pm \Re(z) \geq 0$ near $0$.

\item[(iii)] For $\delta > 0$ fixed, let 
$\Gamma^{\delta}(r,1) \subset \mathcal{D}$ be defined as 
in \eqref{eq2,11}. Assume that there exists a constant $\gamma 
> 0$ such that 
$$
n \big( [r,1] \big) = \mathcal{O}(r^{-\gamma}), \quad r \searrow 0,
$$ 
with $n \big( [r,1] \big)$ growing unboundedly as $r \searrow 0$. 
Then, there exists a positive sequence $(r_{\ell})_\ell$ which 
tends to $0$ such that 
\begin{equation}\label{eqa,9}
\mathcal{N} \big( \overline{\Gamma^{\delta}(r_{\ell},1)} \big) = n 
\big( [r_{\ell},1] \big) \big( 1 + o(1) \big), \quad \ell \rightarrow 
\infty.
\end{equation}

\item[(iv)] If we have
$$
n \big( [r,1] \big) = \Phi(r) \big( 1 + o(1) \big), \quad r 
\searrow 0,
$$ 
with $\phi \big( r(1 \pm \nu) \big) = 
\phi (r) \big( 1 + o(1) + \mathcal{O}(\nu) \big)
$ for any $\nu > 0$ small enough, then
\begin{equation}\label{eqa,10}
\mathcal{N} \big( \overline{\Gamma^{\delta}(r,1)} \big) = 
\Phi(r) \big( 1 + o(1) \big), \qquad r \searrow 0.
\end{equation}
\end{itemize}
\end{lem}


\end{document}